\documentclass[reqno,12pt]{amsart}

\usepackage{psfig, amsmath, amsfonts, amssymb, amsthm, amscd, graphics,epsf,psfrag}

\numberwithin{equation}{section}

\makeatletter
\def\captionfont@{\footnotesize}
\def\captionheadfont@{\scshape}

\long\def\@makecaption#1#2{%
  \vspace{2mm}
  \setbox\@tempboxa\vbox{\color@setgroup
    \advance\hsize-6pc\noindent
    \captionfont@\captionheadfont@#1\@xp\@ifnotempty\@xp
        {\@cdr#2\@nil}{.\captionfont@\upshape\enspace#2}%
    \unskip\kern-6pc\par
    \global\setbox\@ne\lastbox\color@endgroup}%
  \ifhbox\@ne 
    \setbox\@ne\hbox{\unhbox\@ne\unskip\unskip\unpenalty\unkern}%
  \fi
  \ifdim\wd\@tempboxa=\z@ 
    \setbox\@ne\hbox to\columnwidth{\hss\kern-6pc\box\@ne\hss}%
  \else 
    \setbox\@ne\vbox{\unvbox\@tempboxa\parskip\z@skip
        \noindent\unhbox\@ne\advance\hsize-6pc\par}%
\fi
  \ifnum\@tempcnta<64 
    \addvspace\abovecaptionskip
    \moveright 3pc\box\@ne
  \else 
    \moveright 3pc\box\@ne
    \nobreak
    \vskip\belowcaptionskip 
  \fi
\relax
}
\makeatother
\def\writefig#1 #2 #3 {\rlap{\kern #1 truecm
\raise #2 truecm \hbox{#3}}}

\psfigurepath{.:./pictures}

\DeclareMathSymbol{\leqslant}{\mathalpha}{AMSa}{"36} 
\DeclareMathSymbol{\geqslant}{\mathalpha}{AMSa}{"3E} 
\DeclareMathSymbol{\eset}{\mathalpha}{AMSb}{"3F}     
\renewcommand{\leq}{\;\leqslant\;}                   
\renewcommand{\geq}{\;\geqslant\;}                   

\newtheorem{lem}{Lemma}[section]
\newtheorem{pro}{Proposition}[section]
\newtheorem{thm}{Theorem}[section]
\newtheorem{cor}{Corollary}[section]
\newtheorem{defi}{Definition}[section]
\newtheorem{rem}{Remark}[section]

\newcommand{\cA}{\ensuremath{\mathcal A}}

\newcommand{\cC}{\ensuremath{\mathcal C}}

\newcommand{\cE}{\ensuremath{\mathcal E}}
\newcommand{\cF}{\ensuremath{\mathcal F}}

\newcommand{\cH}{\ensuremath{\mathcal H}}

\newcommand{\cL}{\ensuremath{\mathcal L}}

\newcommand{\cP}{\ensuremath{\mathcal P}}

\newcommand{\cR}{\ensuremath{\mathcal R}}
\newcommand{\cS}{\ensuremath{\mathcal S}}
\newcommand{\cT}{\ensuremath{\mathcal T}}

\newcommand{\cZ}{\ensuremath{\mathcal Z}}

\newcommand{\frJ}{\ensuremath{\mathfrak J}}

\newcommand{\frb}{\ensuremath{\mathfrak b}}

\newcommand{\frg}{\ensuremath{\mathfrak g}}

\newcommand{\bbB}{{\ensuremath{\mathbb B}} }

\newcommand{\bbE}{{\ensuremath{\mathbb E}} }

\newcommand{\bbS}{{\ensuremath{\mathbb S}} }

\newcommand{\bbZ}{{\ensuremath{\mathbb Z}} }

\newcommand{\ga}{\alpha}
\newcommand{\gb}{\beta}
\newcommand{\gd}{\delta}
\newcommand{\gep}{\varepsilon}       
\newcommand{\gp}{\varphi}

\newcommand{\gz}{\zeta}

\newcommand{\gD}{\Delta}

\newcommand{\go}{\omega}

\newcommand{\gL}{\Lambda}
\newcommand{\gs}{\sigma}

\newcommand{\Perc}{\Phi}           

\def\1{\ifmmode {1\hskip -3pt \rm{I}} \else {\hbox {$1\hskip -3pt \rm{I}$}}\fi}

\newcommand{\lra}{\leftrightarrow}

\newcommand{\bra}{\langle}
\newcommand{\ket}{\rangle}

\newcommand{\slab}{{\cS}}
\newcommand{\topb}{\partial^{{\rm top}}}
\newcommand{\botb}{\partial^{{\rm bot}}}
\newcommand{\FKm}[2]{\Perc^{#1}_{#2}}

\newcommand{\boxh}{{\bbB}}   
\newcommand{\seed}{\text{seed}}

\newcommand{\north}{\text{North}}
\newcommand{\east}{\text{East}}
\newcommand{\west}{\text{West}}
\newcommand{\south}{\text{South}}

\newcommand{\bond}{p}

\title{Slab percolation for the Ising model}

\author{T. Bodineau}
\address{D{\'e}partement de math{\'e}matiques, Universit{\'e} Paris 7,
case 7012, 2 place Jussieu, Paris 75251, France}
\email{bodineau@gauss.math.jussieu.fr}
\thanks{I wish to thank G. Grimmett, D. Ioffe and R. Kotecky for very 
stimulating discussions and useful comments.}

\date{\today}

\begin{document}

\maketitle

\begin{abstract}
For the FK representation of the Ising model, we prove that 
the slab percolation threshold coincides with the critical temperature    
in any dimension $d \geq 3$.
\end{abstract}

\section{Introduction}

Renormalization arguments are at the core of the description of 
microscopic systems.
In particular, they are necessary close to the critical point
where perturbative techniques no longer apply.
The rigorous implementation of renormalization techniques is model
dependent and in this paper we shall focus on  the 
$q$-Potts model. In this case, 
the critical temperature is characterized by a breaking of symmetry
(spontaneous magnetization) or by the occurrence of percolation
(in the FK representation).
In principle, this characterization should suffice to obtain
further information about the sub-critical or super-critical phases, 
like the control of the susceptibility, the classification of the 
phases and so on.
Unfortunately, some of these issues can only be settled at very high 
or very low temperatures and otherwise they remain open.
In fact, in the intermediate regime of temperatures, the knowledge of
the absence/occurrence of percolation might not be enough to 
implement a mathematical argument.
Nevertheless, progress can be made under stronger assumptions which
sometimes can be proven afterward.\\

The art of statistical mechanics is to propose the right criteria 
from which concrete results can be deduced and 
which should be valid in the whole sub-critical/super-critical regime.
As an example, the Dobrushin and Shlosman strong mixing property
\cite{DS} implies that the sub-critical phase is well behaved
(complete analyticity ...), as well as the dynamics. 
In some instances this property can be checked: for example 
Schonmann and Shlosman \cite{SS} have shown that it is valid
in the whole uniqueness regime of the two-dimensional Ising model 
with nearest neighbor interaction.

In this paper we address another type of hypothesis, namely  
the percolation in slabs.
This concept was introduced in the context of Bernoulli percolation
by Aizenman, Chayes, Chayes, Fr{\"o}hlich, Russo \cite{ACCFR} and
turned out to be a crucial tool to derive many facts about the 
geometry of open clusters (see e.g. \cite{G1}).
Let us be more specific.
If  $p_c$ denotes the critical intensity for the Bernoulli percolation,
then for any $p>p_c$, the origin has a positive probability to be
connected to infinity by an open path in $\bbZ^d$.
On the other hand, one can consider the stronger constraint
that the origin is connected to infinity in a two-dimensional
slice of $\bbZ^d$ for $d \geq 3$.
This happens with positive probability for any $p$ larger than a
critical value $\hat p_c$.
Above $\hat p_c$ several procedures have been developed to analyze
many physical phenomena. Thus,  all that is needed
to generalize these results to the entire 
super-critical regime,  is to prove that $\hat p_c$
and $p_c$ coincide.
This conjecture  was solved in the breakthrough work of 
Grimmett and Marstrand \cite{GM} also
using important ideas introduced by Barsky, Grimmett, Newman 
\cite{BGN}.
The concept of slab percolation was generalized successfully 
to the random cluster measure by Pisztora \cite{pisztora}, who 
proposed an extremely powerful
renormalization scheme under the hypothesis that slab percolation 
occurs.
Furthermore Pisztora conjectured that for FK percolation, 
slab percolation should also be valid in the whole super-critical 
phase.
This coarse graining provides a very good description of the 
super-critical regime and in particular, it is a crucial
tool in the derivation of the Wulff construction \cite{Ce,CePi,Bo,BIV}.\\

The main object of this paper is to prove that for the random cluster
measure associated to the Ising model, the slab percolation
threshold coincides with the critical temperature.
The proof uses heavily the strategy of dynamic renormalization 
introduced  by Barsky, Grimmett, Newman \cite{BGN} in the context of 
Bernoulli percolation.
The starting point of their method was the assumption of percolation
in a half space, from which a coarse graining could be implemented.
For the random cluster measure, the basic characterization of the 
super-critical regime is a positive probability of percolation in any
finite box with wired boundary conditions. The main difficulty 
with implementing this information is created by the dependence on the
boundary conditions.

Our approach relies on the surface tension from which accurate 
estimates on  percolation in finite size volumes can be obtained
uniformly over the boundary conditions.
This requires precise controls on the surface tension which are known
for the Ising model, thanks to inequalities.
The most important property being the positivity of the surface tension
derived by Lebowitz and Pfister \cite{LP} in the whole phase transition
regime. 
The consequences of the surface tension estimates are summarized in
Corollary \ref{cor: connection} from which the coarse graining will be
constructed.
Besides Subsection \ref{subsec: Surface tension estimates}, the 
renormalization procedure developed in the rest of the paper is
valid for general random cluster measures $q \geq 1$.
Further heuristics as well as the scheme of the proof will be
presented in Subsection \ref{subsec: Heuristics}.

\section{Notation and Result}

\subsection{The random cluster measure}
\label{subsec: FK}

Let $\bbE^d$ be the set of bonds, i.e. of pairs  of nearest neighbor vertices
$(i,j)$.
The set $\Omega = \{ 0,1\}^{\bbE^d}$  is the state space 
for the dependent percolation measures.
Given $\go\in\Omega$ and a bond $b=(i,j) \in\bbE^d$, we say  that $b$ is open 
if $\go (b)=1$. Two sites of $\bbZ^d$ are said
to be connected if one can be reached from another via a chain of open bonds.
Thus, each $\go\in\Omega$ splits  $\bbZ^d$ into the disjoint union of maximal
connected components, which are called the open clusters of $\Omega$. Given a
finite subset $B\subset\bbZ^d$  we use $c_B (\go )$ to denote the number of
different open finite 
clusters of $\go$ which have a non-empty intersection
with $B$.

For any $\gL \subset \bbZ^d$ we define the random cluster measure 
on the bond configurations $\omega\in \Omega_\gL  =  \{0
,1\}^{\bbE_\gL}$,  where  $\bbE_\gL$ is the set 
of bonds $b\in\bbE^d$ intersecting $\gL$. 
The boundary conditions are specified by a frozen percolation configuration  
$\pi \in \Omega_\gL^c = \Omega \setminus \Omega_\gL $. 
Using the shortcut  $c^\pi_\gL (\go ) =c_{\gL} (\go\vee \pi )$ for 
the joint configuration $\go \vee\pi\in\bbE^d$, we define the finite 
volume random cluster measure $\FKm{p ,\pi}{\gL}$ on $\Omega_\gL$  with the 
boundary conditions $\pi$ as:
\begin{equation}
\label{FKm}
\FKm{p ,\pi}{\gL}\left(\go \right)~ = ~\frac1{Z^{p ,\pi}_\gL}
\left( \prod_{b\in
\bbE_\gL} \big( 1-p \big)^{1-\go_b} \; p^{\go_b} 
\right)\, q^{c^\pi_\gL (\go)}\,,
\end{equation}
for $q \geq 1$. In this paper, we will sometimes use bond dependent
intensities $p(b)$.
When there is no risk of confusion, we drop the upper-script $p$
and write  $\FKm{\pi}{\gL}$.

\medskip

The random cluster measure associated to the box $\gL_N = \{ -N, \dots, N \}^d$ 
will be denoted by $\FKm{\pi}{N}$.
The measures $\FKm{\pi}{N}$ are FKG partially ordered with respect to
the lexicographical order of the boundary condition $\pi$. Thus, the extremal
ones correspond to the free ($\pi\equiv 0$) and wired ($\pi\equiv 1$) boundary
conditions and are denoted as  $\FKm{{\rm f}}{N}$ and $\FKm{{\rm
w}}{N}$  respectively. 
The  corresponding infinite volume ($N\to\infty$)
limits  $\FKm{{\rm f}}{\,}$ and $\FKm{{\rm w}}{\,}$ always exist.
In the following, the same set of bonds $\bbE_\gL$ will be used 
for random cluster measures in $\gL$ with free or with wired boundary 
conditions.

\medskip

A correspondence between the $q$-Potts model and the random cluster 
measure was established by Fortuin and Kasteleyn \cite{FK}
(see also \cite{ES,G2}).
This representation of the Potts model will be referred as FK
representation.
Of particular interest for us is the Ising model at inverse
temperature $\gb$ which can be related to the
previous model by setting $q =2$ and choosing the bond intensity
$p = 1-\exp(-2\gb)$.
More precisely, the Ising model on $\bbZ^d$ with nearest neighbor 
interaction is defined in terms of spins $\{ \gs_i \}_{ i \in \bbZ^d}$
taking values  $\pm 1$. 
Let $\gs_{\gL_N} \in \{ \pm 1 \}^{\gL_N}$ be the spin configuration 
restricted to $\gL_N$.
The Hamiltonian associated to $\gs_{\gL_N}$ with boundary conditions 
$\gs_{\partial \gL_N}$ is defined by
\begin{eqnarray*}
H( \gs_{\gL_N}  \, | \, \gs_{\partial \gL_N} ) = 
- {1 \over 2} \sum_{ i \sim j \atop i,j \in \gL_N} \gs_i \gs_j
- \sum_{ i \sim j \atop i \in \gL_N, j \in \partial \gL_N} \gs_i \gs_j.
\end{eqnarray*}
The Gibbs measure in $\gL_N$ at inverse temperature $\gb > 0$ is
defined by
\begin{eqnarray*}
\label{Gibbs measure}
\mu_{\gb, \gL_N}^{\gs_{\partial \gL_N}} ( \gs_{\gL_N} ) =
{1 \over Z_{\gL_N}^{\gs_{\partial \gL_N}} }
\exp \big( - \gb  H( \gs_{\gL_N}  \, | \, \gs_{\partial \gL_N} ) \big), 
\end{eqnarray*}
where the partition function $Z_{\gL_N}^{\gs_{\partial \gL_N}}$ is the normalizing 
factor. 
The boundary conditions act as  boundary fields, therefore more general
values of the boundary conditions can be used.

\subsection{The slab percolation threshold}

The phase transition of the random cluster model is characterized by 
the occurrence of percolation above the critical intensity $p_c$
\begin{equation}
\label{eq: transition}
\forall p > p_c, \qquad 
\lim_{N\to\infty}\FKm{p, {\rm w}}{N}\left( 0 \lra
\gL_N^c\right)~=~ \FKm{p, {\rm w}}{\,} \left( 0\lra\infty \right)
> 0.
\end{equation}
If $\gb_c$ denotes the critical inverse temperature of the Ising
model, the FK representation implies that $p_c = 1-\exp(-2\gb_c)$
 for $q =2$.

\medskip

For $d \geq 3$, one may wonder if the stronger property of percolation
in a slab also holds up to the critical value.
For any integers $(L,N)$ we define the slabs of thickness $L$ as
\begin{equation*}
\slab_{L,N} =\{-L, \dots, L \}^{d-2} \times \{-N, \dots, N \}^2 \, , 
\qquad
\slab_L = \{-L, \dots, L \}^{d-2} \times \bbZ^2 \, .
\end{equation*}
A critical value can be associated to any slab thickness $L$
\begin{equation}
\label{eq: slab threshold}
\hat p_c (L) = \inf \left\{ p \geq 0, \qquad 
\liminf_{N} \; \inf_{x \in \slab_{L,N}} \;
\FKm{p, {\rm f}}{\slab_{L,N}} \left( 0\lra x \right) > 0 
\right\} \, .
\end{equation}
As the function $L \to \hat p_c (L)$ is non-increasing, it admits a limit
$\hat p_c$ as $L$ tends to infinity. The critical value
$\hat p_c$ is the {\bf slab percolation threshold} and satisfies
$\hat p_c \geq  p_c$.
The main result of this paper is to prove that both critical
values coincide if $q =2$.

\begin{thm}
\label{thm: slab threshold}
\ 
\begin{itemize}
\item
If $d \geq 3$ and $q=2$, then $\hat p_c =  p_c$.
\item
More generally for $q \geq 1$, if $\Theta_q$ is the set
defined before Corollary \ref{cor: connection} then 
$]\hat p_c,1] \subset \Theta_q$.
\end{itemize}
\end{thm}

To avoid many technicalities, the Theorem has been derived for
random cluster measures with nearest neighbor interactions.
Nevertheless, the proof can be adapted in a straightforward manner
to the case of finite range interactions.
Unbounded range interactions would require a more 
delicate treatment (see e.g. \cite{MS,GH}).

\vskip.4cm

As a by-product of the proof, we get
\begin{cor}
\label{cor: TS}
If $d \geq 3$ and $q=2$, then the surface tension of the Ising model
(see Definition \ref{def: TS}) is equal to 0 at $p_c$.
\end{cor}
This follows from an argument similar to the one used in \cite{BGN}
to prove  the continuity of phase transition for half-space Bernoulli
percolation.

	\subsection{Heuristics and scheme of the proof}
	\label{subsec: Heuristics}

In the super-critical regime, the percolation cluster can be
seen as a backbone on which small open clusters are attached
(the leaves).
When the bond density approaches $p_c$, the backbone becomes 
thinner and the structure of the leaves becomes more chaotic as the
correlation length $\gz(p)$ diverges. 
Nevertheless,  for any $p>p_c$, patterns with similar features 
are repeated on a scale of the order of the correlation length, and
the bond configurations which are distant from each other by at 
least $\gz(p)$ behave essentially independently.
Thus the backbone of the 
percolation cluster has enough space to spread in any slab
of thickness much larger than $\gz(p)$.
This heuristic justifies the slab percolation  conjecture.\\

\begin{figure}[h]
\begin{center}
\leavevmode
\epsfysize = 5 cm
\epsfbox{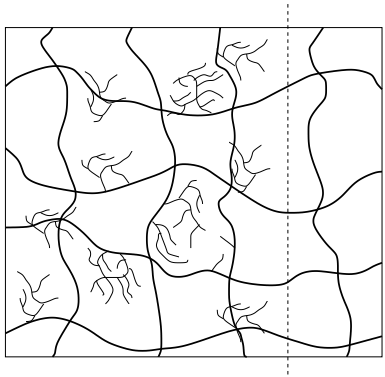}
\end{center}
\caption{A backbone with some leaves is depicted. The dashed line
shows which bonds have to be closed to disconnect the two faces
of the box.}
\label{fig: TS}
\end{figure}

Some information about the backbone structure is encoded in the 
surface tension.
If the surface tension is positive then the probability that one 
face of  a cube of side length $N$   
is not connected to the opposite face decays like $\exp ( - \tau N^{d-1})$.
This means that to disconnect one face from the other a number 
of bonds of the order of $N^{d-1}$ have to be closed
(see Figure \ref{fig: TS}). Thus one deduces 
that the intersection of the backbone with any hyperplan contains
a density of bonds.

The implication of the surface tension on the percolation are
derived in  Subsection \ref{subsec: Surface tension estimates}.
For the Ising model, the positivity of the surface tension 
has been established in the whole super-critical regime \cite{LP}.
Furthermore, in the framework of the Ising model,  several estimates
on the surface tension can be derived uniformly wrt the boundary conditions.
This independence wrt the boundary conditions will be crucial to
implement the renormalization scheme.
With the exception of Subsection \ref{subsec: Surface tension estimates},
the rest of the paper does not rely on the properties of the Ising model.
Thus for any $q \geq 1$, the coarse graining also applies  under
the assumption that $p$ belongs to the range of intensities $\Theta_q$
for which Corollary \ref{cor: connection} holds.
We believe that  the estimates on the surface tension derived for the
Ising model should be valid for any random cluster model if $p>p_c$.\\

This information about the backbone can be used to implement a dynamic
renormalization procedure in the spirit of the one introduced
by  Barsky, Grimmett, Newman \cite{BGN}.
Before going into the details, let us comment on the results
obtained in the case of Bernoulli percolation.
Under the assumption of half-space percolation, Barsky et al.
proposed a renormalization scheme in two steps. First they construct
large blocks, so that with high probability the faces of these 
blocks are interconnected in such a way that further connections
can be launched from any faces of these blocks.
Using the independence between disjoint blocks, they were able 
in a second step of renormalization to create an infinite open cluster
by piling up these blocks.
The simplifying feature of the half-space percolation is to decouple
the different blocks, the drawback being that the
threshold of half-space percolation may not coincide with $p_c$.
Grimmett and Marstrand overcame this problem by using the 
sprinkling technique. 
At the cost of a small increase of bond intensity, they were
able to control the block connections without assuming half-space 
percolation.
When considering the random cluster measure, the dependence on the
boundary conditions adds up and it seems difficult to generalize 
their proof directly from the knowledge that percolation occurs in
any finite box with wired boundary conditions.

The surface tension estimates imply precise controls on the probability
of percolation in a half box uniformly wrt the boundary conditions
(see Subsection \ref{subsec: Half box percolation}).
This will be used to adapt the strategy of Barsky et al. :
the main effort is devoted to performing the first renormalization 
step, i.e. to constructing coarse grained blocks such that with high
probability they satisfy good properties uniformly wrt the boundary 
conditions.
As in \cite{BGN}, the blocks are such that a small region on the
underside is connected to seeds lying in every facet of the 
block.
Nevertheless, the lack of independence wrt the boundary conditions
has prevented us from directly applying the ideas of  \cite{BGN}
and several detours are necessary.
This task is performed in Section \ref{sec:  first renormalization step}.
Once this is done the coarse grained blocks can be piled up as if
they were independent and the second renormalization step applies 
as in \cite{BGN}.
For the sake of completeness, the main ideas of the geometric
construction, including the steering and branching  rules,
are recalled in Section \ref{sec: second renormalization level}.

\vskip.2cm

Finally, we would like to comment on two related works based on 
very different strategies.

The convergence of the critical temperatures of three-dimensional
slabs has been derived by Aizenman in the context of Ising 
model \cite{A1,A2}.
The framework is different from ours and in particular the slabs 
are of the form $\bbZ^3 \times \{-L,L\}$ and the boundary conditions 
are periodic.
As the proof relies on the random current representation of the Ising
model and on reflection positivity, it seems difficult to relate 
Aizenman's results to the slab percolation threshold defined in the
FK setting \eqref{eq: slab threshold} and therefore to Pisztora's
coarse graining.

For very large $q$, the FK measure can be analyzed by means of
Pirogov Sinai theory (see e.g. \cite{BKM,KLMR}). 
At $p_c$ the transition is
first order and above $p_c$, one can show that the only stable phase 
is the one with high density of open bonds.
Thus for large $q$, the slab percolation threshold becomes a trivial
matter and renormalization techniques do not seem necessary to
describe the super-critical phase\footnote{Private communication by 
R. Kotecky.}.

\vskip1cm

\section{Crossing clusters}

        \subsection{Surface tension estimates}
	\label{subsec: Surface tension estimates}

In this Subsection we study the dependence on the boundary conditions of 
the surface tension.

Let us first recall the definition of the surface tension along the 
coordinate axis $\vec{e}_d$.
Let $\gd$ be a positive parameter (typically chosen very small) and $N,L$
be two integers such that $N \gg L$.
In the following these parameters are chosen of the form
$N = 2^n, L = 2^\ell, \gd = 2^{-p}$, where $n,\ell,p$ are integers.
We consider the increasing family of rectangles 
\begin{eqnarray*}
\cR^L (N,\gd) = \{-N,\dots,N\}^{d-1} \times \{ - \gd N - L,\dots,\gd N + L \} \, ,
\end{eqnarray*}
and simply write $\cR (N,\gd)$ when $L =0$. 
Finally, we introduce the set $\frJ(N,\gd)$ of bond configurations 
for which there is no connection from the
top face to the bottom face of $\cR (N,\gd)$ (i.e. the two faces orthogonal 
to $\vec{e}_d$)
\begin{eqnarray}
\frJ(N,\gd) = \left\{
\botb \cR (N,\gd) \ \not \lra \  \topb \cR (N,\gd) \right\} \, .
\end{eqnarray}

\begin{defi}
\label{def: TS}
The surface tension in the direction $\vec{e}_d$ is defined by
\begin{eqnarray*}
\tau_p  = \lim_{\gd \to 0} \;
\lim_{N \to \infty} - \frac{1}{N^{d-1}}
\log \FKm{p,{\rm w}}{\cR (N,\gd)} \Big( \frJ(N,\gd) \Big) \, .
\end{eqnarray*}

For $q =2 $ and $p = 1 - \exp( -2\gb)$, the Ising counterpart is
\begin{eqnarray*}
\tau_p  = \lim_{\gd \to 0} \;
\lim_{N \to \infty} - \frac{1}{N^{d-1}}
\log \frac{Z^\pm_{\cR (N,\gd)}}{Z^+_{\cR (N,\gd)}} \, ,
\end{eqnarray*}
where $Z^\pm_{\cR (N,\gd)}$ denotes the partition function with 
mixed boundary conditions.
\end{defi}

The convergence of the thermodynamic limit has been derived in
\cite{MMR}.
Notice also that the lateral sides have a vanishing perimeter 
and therefore they have no influence on the value of the 
surface tension.

\vskip.5cm

For our purposes,
it will be useful to consider the equivalent formulation 
\begin{eqnarray}
\label{def: surface tension}
\tau_p
=\lim_{\gd \to 0} \;
\lim_{N \to \infty} - \frac{1}{N^{d-1}}
\log \FKm{p,{\rm w}}{\cR^L (N,\gd)} \Big( \frJ(N,\gd) \Big)
\, .
\end{eqnarray}
The dependence on $L$ can be estimated by using FKG inequality and
the fact that  $\frJ(N,\gd)$ is a non increasing event
\begin{eqnarray*}
\FKm{{\rm w}}{\cR (N,\gd)} \Big( \frJ(N,\gd) \Big) 
\leq
\FKm{{\rm w}}{\cR^L (N,\gd)} \Big( \frJ(N,\gd) \Big)
\leq
\FKm{{\rm w}}{\cR^L (N,\gd)} \Big( \left\{
\topb \cR^L (N,\gd) \ \not \lra \ \botb \cR^L (N,\gd) \right\} 
\Big)
\, .
\end{eqnarray*}
As the LHS and the RHS (properly renormalized) converge to the surface tension,
the identity \eqref{def: surface tension} holds.

\vskip1cm

The following result enables us to compare the surface tension for different 
boundary conditions. The parameter $L$ should be taken large enough in order to 
screen the boundary conditions.
\begin{thm}
\label{thm: equality ST} 
For $q =2$, then for any $p \in [0,1]$,
\begin{eqnarray}
\label{eq: equality ST}
\left|
\log \FKm{p, {\rm w}}{\cR^L (N,\gd)} \Big( \frJ(N,\gd) \Big)
- 
\log \FKm{p, {\rm f}}{\cR^L (N,\gd)} \Big( \frJ(N,\gd) \Big)
\right|
\leq \left( \gep_L + c_p \gd + c_d \frac{L}{N} \right) N^{d-1} \, ,
\end{eqnarray}
where $\gep_L$ vanishes as $L$ tends to infinity.
\end{thm}

The specificity of the Ising model will be used only for the
derivation of \eqref{eq: connection free}.
Thus for any $q \geq 1$, we introduce the set of bond intensities
$\Theta_q$ for which \eqref{eq: connection free} holds.
The rest of the paper will apply for any random cluster measure
with $p \in \Theta_q$.

We can now state the main result of this  Subsection.
\begin{cor}
\label{cor: connection}
Let $q =2$, then for any $p > p_c$, there exists 
$\gd = \gd(p)$ and $L$ large such that for any $N$ large enough
\begin{eqnarray}
\label{eq: connection free}
\FKm{p, {\rm f}}{\cR^L (N,\gd)} \Big( 
\left\{\botb \cR (N,\gd) \ \leftrightarrow \ \topb \cR (N,\gd) \right\}
\Big) 
\geq 1 - \exp( - C_p  N^{d-1} ) \, ,
\end{eqnarray}
where $C_p$ is a positive constant.
\end{cor}
This follows from Theorem \ref{thm: equality ST} and the positivity 
of surface tension
which was proven for the Ising model by Lebowitz and Pfister 
\cite{LP} when $\gb$ is larger than $\gb_c$.

\vskip1cm

\noindent
{\it Proof of Theorem \ref{thm: equality ST}}.

The estimate \eqref{eq: equality ST} will be obtained by interpolating
the boundary conditions wired and free.
Define the set of boundary bonds at the top and bottom faces of
$\cR^L (N,\gd)$ as
\begin{eqnarray*}
\Xi = \left\{ (i,i+\vec{e}_d), (j,j-\vec{e}_d) \in 
\bbE_{\cR^L (N,\gd)}, 
\qquad
i,j \in \cR^L (N,\gd); \ i_d = - j_d = L + \gd N 
\right\} \, .
\end{eqnarray*}
Let $\FKm{s,{\rm w}}{\cR^L (N,\gd)}$ be the wired FK measure
for which the bonds in $\Xi$ have intensity $s$
instead  of $\bond$.
The parameter $s$ acts as a boundary magnetic field (on the faces orthogonal 
to $\vec{e}_d$) and interpolates between $\FKm{{\rm w}}{\cR^L
(N,\gd)}$ (for $s= \bond$) and 
$\FKm{{\rm f,w}}{\cR^L (N,\gd)}$ (for $s=0$).
The latter measure is the FK measure with free boundary conditions
on the top and bottom faces of $\cR^L (N,\gd)$ and
wired otherwise.

\begin{rem}
In Subsection \ref{subsec: FK}, the FK measure $\FKm{{\rm f,w}}{\cR^L
(N,\gd)}$ was  defined  on the set of bonds intersecting $\cR^L (N,\gd)$.
When $s = 0$, the measure $\FKm{s=0,{\rm w}}{\cR^L (N,\gd)}$ is 
equal to the free measure conditionally to the
fact that the bonds in $\Xi$ are closed.
Nevertheless this has no impact on the probability of events
which are not supported
by $\Xi$ and thus the probability of $\frJ(N,\gd)$ is the same under 
$\FKm{{\rm f,w}}{\cR^L (N,\gd)}$ or $\FKm{s=0,{\rm w}}{\cR^L (N,\gd)}$.
\end{rem}

\begin{figure}
\begin{center}
\leavevmode
\epsfysize = 5 cm
\psfrag{b}[Br]{$\cR (N,\gd)$}
\psfrag{a}[Br]{$\cR^L (N,\gd) \setminus \cR (N,\gd)$}
\psfrag{c}{$\Xi$}
\psfrag{d}{$L$}
\psfrag{e}{$2 \gd N$}
\epsfbox{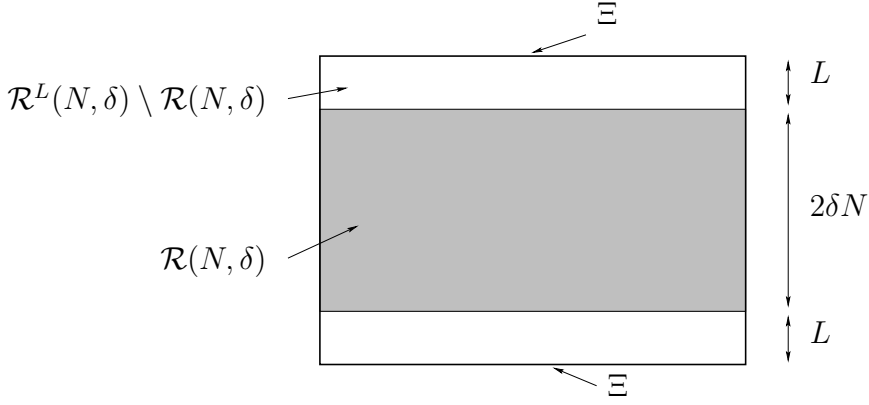}
\end{center}
\caption{The event $\frJ(N,\gd)$ is supported by the shaded region
which is decoupled from $\Xi$.}
\label{fig: TS}
\end{figure}

\begin{eqnarray}
\label{eq: 1.1}
&&\log \FKm{{\rm w}}{\cR^L (N,\gd)} \Big( \frJ(N,\gd) \Big)
-
\log \FKm{{\rm f,w}}{\cR^L (N,\gd)} \Big( \frJ(N,\gd) \Big)\\
&& \qquad \qquad =
\sum_{b \in \Xi} \int_0^{\bond}  \frac{ds}{s(1- s)} \,
\left(
\frac{\FKm{s,{\rm w}}{\cR^L (N,\gd)} \Big( \frJ(N,\gd) \, \go_b 
\Big) }{\FKm{s,{\rm w}}{\cR^L (N,\gd)}\Big( \frJ(N,\gd) \Big)}
-
\FKm{s,{\rm w}}{\cR^L (N,\gd)} (  \go_b )
\right)
\, . \nonumber
\end{eqnarray}

At this stage, a mixing property is required.  We  introduce
now the FK measures on the box 
$\bbB_K = \{ -K,\dots, K\}^{d-1} \times \{ 0, \dots, K\}$ 
with intensity $s$ for the bonds below the bottom face of $\bbB_K$ 
(i.e. the bonds in $\{ (i,i-\vec{e}_d), \quad i_d = 0, i \in  \bbB_K \}$) 
and wired boundary conditions at the bottom face.
When the boundary conditions on the remaining faces are free the
corresponding  FK measure is denoted by $\FKm{s,{\rm f}}{\bbB_K}$
 and $\FKm{s,{\rm w}}{\bbB_K}$ for wired.
\begin{pro}
\label{prop: mixing}
Let $\frb_0$ be the bond $(0,-\vec{e}_d)$.
For any $\bond \in [0,1]$ and $s \in (0, \bond]$,
\begin{eqnarray}
\label{eq: mixing}
\left| 
\FKm{s,{\rm w}}{\bbB_K} ( \go_{\frb_0} )
-
\FKm{s,{\rm f}}{\bbB_K} ( \go_{\frb_0} ) 
\right|
\leq \gep_K (s) \, ,
\end{eqnarray}
where $\gep_K (s)$ vanishes as $K$ diverges.
\end{pro}

We postpone the proof of this Lemma and first estimate  \eqref{eq: 1.1}.
As $\frJ(N,\gd)$ is supported by $\cR (N,\gd)$, it is decoupled from
the bonds in $\Xi$.
Thus using FKG property one can reduce the estimates
to the set $\cR^L (N,\gd) \setminus \cR (N,\gd)$
\begin{eqnarray*}
&&
\left| \log \FKm{{\rm w}}{\cR^L (N,\gd)} \Big( \frJ(N,\gd) \Big)
-
\log \FKm{{\rm f,w}}{\cR^L (N,\gd)} \Big( \frJ(N,\gd) \Big)
\right| \\
&& \qquad \qquad \leq
\sum_{b \in \Xi} \int_0^{\bond} \frac{ds}{s(1-  s)} \,
\left( 
\FKm{s,{\rm w}}{\cR^L (N,\gd) \setminus \cR (N,\gd)} ( \go_b )
- 
\FKm{s,{\rm f}}{\cR^L (N,\gd) \setminus \cR (N,\gd)} ( \go_b )
\right)  \, , \\
&& \qquad \qquad \leq
2 N^{d-1} \int_0^{\bond}  \frac{ds}{s(1-  s)} \,
\left( 
\FKm{s,{\rm w}}{\bbB_L} ( \go_{\frb_0} )
- 
\FKm{s,{\rm f}}{\bbB_L} ( \go_{\frb_0} )
\right) + c_d L N^{d-2} \, .
\end{eqnarray*} 
The final bound has been obtained by applying again FKG inequality to
reduce further the domain of the FK measure. The rest is an upper
bound of the expectation of the terms lying at a distance smaller than
$L$ from the lateral sides of  $\cR^L (N,\gd)$.\\

We remark that for any $s \in (0,\bond]$
\begin{eqnarray*}
&&  \frac{1}{s(1-  s)} \,
\left( 
\FKm{s,{\rm w}}{\bbB_L} ( \go_{\frb_0} )
- 
\FKm{s,{\rm f}}{\bbB_L} ( \go_{\frb_0} )
\right)
\leq
\frac{1}{s(1-  s)} \,
\left( 
\FKm{s,{\rm w}}{\{ \frb_0 \}} ( \go_{\frb_0} )
- 
\FKm{s,{\rm f}}{ \{ \frb_0 \} } ( \go_{\frb_0} )
\right) \, ,\\
&& \qquad  \leq   \frac{1}{s(1-  s)} \,
\left( s - \frac{s}{s + (1-s) q} \right) 
\leq \frac{1}{1-  \bond} \, .
\end{eqnarray*}

Thus the dominated convergence theorem and Proposition \ref{prop: mixing}
imply that there exists a sequence $\gep_L$ vanishing to 0 as $L$
diverges such that
\begin{eqnarray*}
\left| \log \FKm{{\rm w}}{\cR^L (N,\gd)} \Big( \frJ(N,\gd) \Big)
-
\log \FKm{{\rm f,w}}{\cR^L (N,\gd)} \Big( \frJ(N,\gd) \Big)
\right|
\leq  2 N^{d-1} \gep_L + c_d L N^{d-2} \, . 
\end{eqnarray*}

It remains to compare the FK measures $\FKm{{\rm f,w}}{\cR^L (N,\gd)}$
and $\FKm{{\rm f}}{\cR^L (N,\gd)}$ by modifying the boundary
conditions on the lateral sides. This can be 
achieved at a finite cost proportional to the perimeter  of the sides
parallel to $\vec{e}_d$.
The corresponding error is bounded by the term $\gd N^{d-1}$
in the LHS of \eqref{eq: equality ST}.
This completes Theorem \ref{thm: equality ST}.  \qed

\vskip1cm

\noindent
{\it Proof of  Proposition \ref{prop: mixing}}.

The strategy is to use the Gibbssian counterpart of \eqref{eq: mixing}
for which a mixing property is known.\\

\noindent
{\bf Step 1.}

We denote by $\partial \bbB_K$ the set of boundary vertices of $\bbB_K$
and define the set of vertices lying in the bottom face of 
$\partial \bbB_K$ as
\begin{eqnarray*}
\partial^{s} \bbB_K &=& \left\{ i, \qquad i_d = -1, \ 
-K \leq i_\ell \leq K \ \ \text{if} \ \ \ell < d \right\} \,  .
\end{eqnarray*} 

The bond $\frb_0$ links the site 0 to the site
$\frg = (0,\dots,0,-1)$.
The event that 0 is connected to the boundary without using the
bond $\frb_0$ is in the wired and free case
\begin{eqnarray*}
\cC^{{\rm w}} =
\left\{ 0 \lra \partial \bbB_K \setminus \{ \frg \} \right\}
\quad \text{and} \quad
\cC^{{\rm f}} = 
\left\{ 0 \lra \partial^{s} \bbB_K \setminus \{ \frg \}
\right\} \, .
\end{eqnarray*}

\medskip

Let $h = h(s) = - \frac{1}{2} \log (1-s)$ be the positive magnetic
field associated to the bond intensity $s$.
\begin{eqnarray*}
\mu^{h,+}_{\bbB_K} ( \gs_0) =
\FKm{s,{\rm w}}{\bbB_K} ( 0 \lra \partial \bbB_K )
= \FKm{s,{\rm w}}{\bbB_K} \big( 1_{\cC^{{\rm w}}}  \big)
+ \FKm{s,{\rm w}}{\bbB_K} \big( ( 1- 1_{\cC^{{\rm w}}}) \,
\go_{\frb_0} 
\big) \, .
\end{eqnarray*} 
Conditioning outside $\{\frb_0\}$, we have for any boundary
condition $\pi$ which does not belong to $\cC^{{\rm w}}$
(see e.g. equation (3.10) in \cite{G2}) 
\begin{eqnarray*}
\FKm{s,\pi}{ \{ 0 \} } \big( \go_{\frb_0} \big)
= \frac{s}{s + (1-s) q} \, .
\end{eqnarray*} 
Thus 
\begin{eqnarray*}
\mu^{h,+}_{\bbB_K} ( \gs_0) =
\FKm{s,{\rm w}}{\bbB_K} \big( 1_{\cC^{{\rm w}}}  \big)
+ 
\left(
1 - \FKm{s,{\rm w}}{\bbB_K} \big( 1_{\cC^{{\rm w}}} \big) 
\right)
\frac{s}{s + (1-s) q}\, .
\end{eqnarray*}

In the same way, with free boundary conditions (i.e. with zero
magnetic field)
\begin{eqnarray*}
\mu^{h,0}_{\bbB_K} ( \gs_0) =
\FKm{s,{\rm f}}{\bbB_K} \big( 1_{\cC^{{\rm f}}}  \big)
+ 
\left(
1 - \FKm{s,{\rm f}}{\bbB_K} \big( 1_{\cC^{{\rm f}}} \big) 
\right)
\frac{s}{s + (1-s) q}\, .
\end{eqnarray*}
This leads to
\begin{eqnarray*}
\FKm{s,{\rm w}}{\bbB_K} \big( 1_{\cC^{{\rm w}}}  \big)
- \FKm{s,{\rm f}}{\bbB_K} \big( 1_{\cC^{{\rm f}}} \big)
= 
\frac{s + (1-s) q}{(1-s) q}
\left(
\mu^{h,+}_{\bbB_K} ( \gs_0)
-
\mu^{h,0}_{\bbB_K} ( \gs_0) 
\right)
\, .
\end{eqnarray*}

\medskip

Thus we have established the following equivalence
\begin{lem}
\label{lem: equivalence}
For $h = - \frac{1}{2} \log (1-s)$, if the difference 
$\big( \mu^{h,+}_{\bbB_K} (\gs_0) - \mu^{h,0}_{\bbB_K} (\gs_0) \big)$ 
vanishes to 0 as $K$ diverges then there exists a sequence 
$\gep_K(s)$ such that
\begin{eqnarray*}
\left| 
\FKm{s,{\rm w}}{\bbB_K} \big( \cC^{{\rm w}}  \big)
- \FKm{s,{\rm f}}{\bbB_K} \big( \cC^{{\rm f}}  \big)
\right|
\leq \gep_K (s) \, ,
\end{eqnarray*}
and $\lim_K \ \gep_K (s) = 0$.
\end{lem}

\vskip.2cm

\noindent
{\bf Step 2.}

Let us prove that for any $h>0$
\begin{eqnarray}
\label{eq: Ising mixing} 
\lim_{K \to \infty}
\mu^{h,+}_{\bbB_K} (\gs_0) - \mu^{h,0}_{\bbB_K} (\gs_0)
= 0 \, .
\end{eqnarray}
The above thermodynamic limits are taken in the upper half space.
We denote by $\mu^{h,+} (\gs_0)$ and $\mu^{h,0} (\gs_0)$
the corresponding limits.

From \cite{MMP}, we know that the functions 
$\mu^{h,+}_{\bbB_K} (\gs_0)$, $\mu^{h,0}_{\bbB_K} (\gs_0)$
as well as $\mu^{h,+} (\gs_0)$ and $\mu^{h,0} (\gs_0)$
are analytic in the complex domain $\{ \text{Re}(h) > |\text{Im}(h)| \}$.
Therefore the functions below enjoy the same properties
\begin{eqnarray*}
F_K (h) = \mu^{h,+}_{\bbB_K} (\gs_0) - \mu^{h,0}_{\bbB_K} (\gs_0)
\quad \text{and} \quad
F (h) = \mu^{h,+} (\gs_0) - \mu^{h,0} (\gs_0)
\, .
\end{eqnarray*}

When $h>1$, then $F(h)= 0$. This follows from an estimate 
derived by  Fr{\"o}hlich and Pfister 
(see (2.13) and (2.21) of \cite{FP})
\begin{eqnarray*}
\forall h>1, \qquad 
\lim_{K \to \infty} \;
\mu^{h,+}_{\bbB_K} (\gs_0) 
- \mu^{h,-}_{\bbB_K} (\gs_0)
= 0 \, .
\end{eqnarray*}
From the FKG inequality, the LHS dominates $F_K (h)$ and
one obtains that $F(h) =0$ for any $h >1$.
As $F$ is analytic it must be constant for any $h>0$.
Thus \eqref{eq: Ising mixing} holds.

\vskip.5cm

\noindent
{\bf Step 3.}

Finally, it remains to control the expectation of 
$\go_{\frb_0}$
\begin{eqnarray*}
\FKm{s,{\rm w}}{\bbB_K} (\go_{\frb_0}) = 
\FKm{s,{\rm w}}{\bbB_K} \big( 1_{\cC^{{\rm w}}} \go_{\frb_0}  \big)
+ \FKm{s,{\rm w}}{\bbB_K} \big( ( 1 - 1_{\cC^{{\rm w}}}) \go_{\frb_0}
 \big)
\end{eqnarray*}
Conditionally to $\cC^{{\rm w}}$ or its complement, the distribution 
of $\go_{\frb_0}$ is determined thus
\begin{eqnarray*}
\FKm{s,{\rm w}}{\bbB_K} (\go_{\frb_0}) =  s
\FKm{s,{\rm w}}{\bbB_K} \big( \cC^{{\rm w}}  \big)
+ \left( 1 - \FKm{s,{\rm w}}{\bbB_K} \big( \cC^{{\rm w}} )
\right) \frac{s}{s + (1 -s) q} \, .
\end{eqnarray*}
A similar relation holds in the case of free boundary conditions,
this leads to
\begin{eqnarray*}
\FKm{s,{\rm w}}{\bbB_K} (\go_{\frb_0}) 
- \FKm{s,{\rm f}}{\bbB_K} (\go_{\frb_0}) 
=  s \left( 1 - \frac{1}{s + (1 -s) q} \right)
\left(
\FKm{s,{\rm w}}{\bbB_K} \big( \cC^{{\rm w}} )
 - \FKm{s,{\rm f}}{\bbB_K} \big( \cC^{{\rm f}} )
\right) 
\, .
\end{eqnarray*}

\medskip

Combining Lemma \ref{lem: equivalence} and \eqref{eq: Ising mixing} 
we conclude Proposition \ref{prop: mixing}.    \qed

\vskip1cm

        \subsection{Crossing clusters}

Corollary \ref{cor: connection} can be improved by 
showing that a connection from the top to the bottom of $\cR^L (N,\gd)$
occurs with high probability.

\begin{pro}
\label{pro: improved connection free}
For any $\bond \in \Theta_q$, there exists 
$\gd$ and $L$ such that for any $N$ large enough
\begin{eqnarray}
\label{eq: improved connection free}
\FKm{p,{\rm f}}{\cR^L (N,\gd)} \Big( 
\left\{\botb \cR^L (N,\gd) \ \leftrightarrow \ \topb \cR^L (N,\gd) \right\}
\Big) 
\geq 1 - \exp( - C  N^{d-1} ) \, ,
\end{eqnarray}
where $C = C(p,L,\gd)$ is a positive constant.
\end{pro}

\noindent
{\it Proof.}

The parameters $\gd$ and $L$ are fixed according to Corollary 
\ref{cor: connection}.
We first analyze the connections from $\topb \cR (N,\gd)$
to the bottom face of $\cR^L (N,\gd)$
and show that there exists $C>0$ such that
\begin{eqnarray}
\label{eq: improved connect 0}
\FKm{{\rm f}}{\cR^L (N,\gd)} \Big( 
\left\{\botb \cR^L (N,\gd) \ \leftrightarrow \ \topb \cR (N,\gd) \right\}
\Big) 
\geq 1 - \exp( - C  N^{d-1} ) \, .
\end{eqnarray}

Given $\ga>0$, we introduce $\cA^\ga_N$ the set 
of bond configurations for which there are at least $\ga N^{d-1}$ 
sites at the level $\{ x_d = - \gd N + 1 \}$ connected to 
$\topb \cR (N,\gd)$ by open paths contained within
$\{ x_d \geq - \gd N + 1 \}$.
There exists $\ga>0$ and $C_1 >0$ such that 
\begin{eqnarray}
\label{eq: connection sites 0}
\FKm{{\rm f}}{\cR^L (N,\gd)} \Big( \cA^\ga_N \Big) 
\geq 1 - \exp( - C_1  N^{d-1} ) \, .
\end{eqnarray}
To see this, we first notice that
the set $(\cA^\ga_N)^c$ is supported by the bonds lying in 
$\{ x_d \geq  - \gd N + 1 \}$.
Given any configuration $\go$ in $(\cA^\ga_N)^c$, it is enough to
close at most $\ga N^{d-1}$ bonds in order to cut all the connections
from $\topb \cR (N,\gd)$ to $\botb \cR (N,\gd)$. 
Therefore, there exists a constant $C_2$ such that
\begin{eqnarray}
\FKm{{\rm f}}{\cR^L (N,\gd)} \Big( (\cA^\ga_N)^c \Big) 
\leq \exp( C_2 \ga   N^{d-1} ) \; \FKm{{\rm f}}{\cR^L (N,\gd)} 
\big( \frJ (N,\gd) \big) \, .
\end{eqnarray}
Using the inequality \eqref{eq: connection free}, one can choose 
$\ga$ small enough such that \eqref{eq: connection sites 0} holds.

\vskip.5cm

The next step is to extend the connections below the level
$\{ x_d = - \gd N + 1 \}$.
For any bond configuration $\pi$ in $\cA^\ga_N$,
we denote by $\{ x^{(i)} \}_i$ the first $\ga N^{d-1}$ sites 
(wrt the lexicographic order) in
$\{ x_d = - \gd N + 1 \}$ which are connected to $\topb \cR (N,\gd)$.
Let ${\bf C} (x^{(i)})$ be the event that $x^{(i)}$ is
connected to $\botb \cR^L (N,\gd)$ by a straight vertical path, i.e. 
using only the edges in
$$
\cP (x^{(i)}) 
= \left( \bra x^{(i)} - (k-1) \vec{e}_d, \; x^{(i)}- k \vec{e}_d \ket
\right)_{k \leq L}
 \, .
$$
Uniformly over the boundary conditions $\go$ outside 
$\cP (x^{(i)})$, we have
\begin{eqnarray}
\label{eq: path 1}
\FKm{\go}{\cP (x^{(i)})} 
\Big( {\bf C} (x^{(i)}) \Big) 
\geq
\FKm{{\rm f}}{\cP (x^{(i)})} 
\Big( {\bf C} (x^{(i)}) \Big) 
\geq c_L >0 \, .
\end{eqnarray}

To simplify the notation, we set $\gD = \cR^L (N,\gd) \cap
\{ x_d \leq -\gd N\}$. 
For any bond configuration $\pi \in \cA^\ga_N$ supported by
$\{ x_d \geq -\gd N +1\}$, we get
\begin{eqnarray*}
\FKm{\pi}{\gD} \Big( 
\left\{\botb \cR^L (N,\gd) \ \lra \ \topb \cR (N,\gd) \right\}
\Big) 
&\geq& 
\FKm{\pi}{\gD} \Big( \bigcup_i  {\bf C} (x^{(i)}) \Big)\\ 
&\geq& 
1 - \FKm{\pi}{\gD} 
\left( \bigcap_i  \left( {\bf C} (x^{(i)}) \right)^c \right) 
\, .
\end{eqnarray*}
The supports of the events ${\bf C} (x^{(i)})$ are disjoint.
Using repeated conditionings and the lower bound \eqref{eq: path 1},
we get uniformly over $\pi$ in $\cA^\ga_N$
\begin{eqnarray*}
\FKm{\pi}{\gD} \Big( 
\left\{\botb \cR^L (N,\gd) \ \lra \ \topb \cR (N,\gd) \right\}
\Big) 
\geq 1 - ( 1 -  c_L)^{\ga N^{d-1}} \geq 
1 - \exp( - C_L \; \ga N^{d-1} ) 
\, , 
\end{eqnarray*}
for some $C_L  >0$.

Finally, combining the previous estimate and 
\eqref{eq: improved connect 0}, there is $C >0$ such that
\begin{eqnarray*}
&& \FKm{{\rm f}}{\cR^L (N,\gd)} \Big( 
\left\{\botb \cR^L (N,\gd) \ \lra \ \topb \cR (N,\gd) \right\}
\Big) \\
&& \qquad \geq
\FKm{{\rm f}}{\cR^L (N,\gd)} \left[ 
\FKm{\pi}{\gD} \Big( 
\left\{\botb \cR^L (N,\gd) \ \lra \ \topb \cR (N,\gd) \right\}
\Big)
\; 1_{\cA^\ga_N} (\pi) \right] \\
&& \qquad \geq 
\FKm{{\rm f}}{\cR^L (N,\gd)} \Big( \cA^\ga_N \Big)
\left(  1 - \exp( - C_L \; \ga N^{d-1} )  \right)
\geq  1 - \exp( - C \, \ga N^{d-1} )
\, .
\end{eqnarray*}
Thus \eqref{eq: improved connect 0} is satisfied. 
The Theorem can be completed by following the same strategy to
extend the connections to the top face of $\cR^L (N,\gd)$.

        \qed

        \subsection{Half box percolation}
	\label{subsec: Half box percolation}

The surface tension estimates imply not only the occurrence of half 
space percolation, but also uniform bounds wrt the boundary conditions
for the  corresponding finite size events. 
At this stage, it will be easier to consider sets in the half space 
$\{ x_d \geq 0\}$ instead of the rectangles $\cR^L (N,\gd)$.
We define below the corresponding notation which will be used throughout
the paper.

\vskip.3cm

\begin{defi}
\ 

\begin{enumerate}
\item
For any integers $(\ell,h)$, a block is defined as
\begin{eqnarray*}
\label{eq: block}
\boxh(\ell,h) = \{ -\ell,\dots ,\ell \}^{d-1} \times \{ 0, \dots, h  \} \, .
\end{eqnarray*}
The interior of $\boxh(\ell,h)$ is 
\begin{eqnarray*}
\boxh^* (\ell,h) = \big\{ -\ell+1, \dots, \ell-1 \big\}^{d-1}
 \times \big\{ 1, \dots, h-1 \} \, .
\end{eqnarray*}
\item
The top face of the block $\boxh(\ell,h)$ will be denoted by 
\begin{eqnarray*}
\label{eq: brick top}
T(\ell,h) = \{ -\ell,\dots ,\ell \}^{d-1} \times \{ h \} \, .
\end{eqnarray*}
The top face is split into $2^{d-1}$ sub-regions
\begin{eqnarray*}
\label{eq: brick tops}
T_1 (\ell,h) &=& \{ 0,\dots ,\ell \}^{d-1} \times  \{  h \} \, ,  \\
T_2 (\ell,h) &=& \{ 0,\dots ,\ell \}^{d-2} \times  
\{ -\ell,\dots ,0 \} \times \{ h \} \, , \nonumber \\
&& \dots \nonumber \\
T_{2^{d-1}} (\ell,h) &=& \{ -\ell,\dots ,0 \}^{d-1} \times  \{ h \} \, . \nonumber
\end{eqnarray*}
\item
Denote by $S (\ell,h)$ the union of the sides
of $\boxh(\ell,h)$. It is divided into $2 (d-1) 2^{d-2}$ 
sub-regions
\begin{eqnarray*}
\label{eq: brick sides}
S_1 (\ell,h) &=& \{\ell\} \times \{ 0,\dots ,\ell \}^{d-2} \times  \{0, \dots,
h \} \, ,  \\
S_2 (\ell,h) &=& \{\ell\} \times \{ 0,\dots ,\ell \}^{d-3} \times  
\{ -\ell,\dots ,0 \} \times \{0, \dots, h \} \, , \nonumber \\
&& \dots \nonumber \\
S_{2 (d-1) 2^{d-2}} (\ell,h) &=& \{ -\ell,\dots ,0 \}^{d-2} \times  \{\ell\} \times 
\{0, \dots, h \} \, . \nonumber
\end{eqnarray*}
\end{enumerate}
\end{defi}

\vskip.5cm

For any $i \leq d$ and $x \in \bbZ^d$, let $b^i_K (x)$ be the 
$(d-1)$-dimensional hypercube centered around $x$ and orthogonal to 
$\vec{e}_i$
\begin{eqnarray*}
b^i_K (x) = \big\{ y \in \bbZ^d, \quad 
y_i = x_i, \ \forall j \ne i, \ | y_j - x_j | \leq K 
\big\} \, .
\end{eqnarray*}
In the following we will be interested to the connections from
$b^i_K (0)$ to the faces of $\boxh (\ell,h)$ by open paths strictly
within $\boxh (\ell,h)$, i.e. paths lying in $\bbE_{\boxh^* (\ell,h)}$.
Recall that $\bbE_{\boxh^* (\ell,h)}$ denotes the set of bonds
intersecting $\boxh^* (\ell,h)$ (see Subsection \ref{subsec: FK}).

\begin{pro}
\label{pro: seed 0}
For any $\bond \in \Theta_q$, there exists $\gd>0$ and $L$
such that the following holds uniformly in $N$ and 
$K \leq N$
\begin{eqnarray}
\label{eq: connection property}
\FKm{p,{\rm f}}{\boxh^* (N,H_N)} \Big( 
\left\{ b^d_K (0)\ \leftrightarrow \ T (N,H_N) \right\}
\Big) 
\geq 1 - \exp( - C K^{d-1} ) \, ,
\end{eqnarray}
where $H_N = \gd N + L$ and $C$ is a positive constant.
\end{pro}

\noindent
{\it Proof.} 

Let $B_N = \big\{ -N, \dots, N \big\}^{d-1} \times \big\{ 0 \big\}$
be the bottom face of $\bbB^* (N,H)$.
According to Proposition \ref{pro: improved connection free}
the following holds for appropriate choices of $\gd >0$ and $L$
\begin{eqnarray}
\label{eq: intermediaire 0}
\FKm{{\rm f}}{\boxh^* (N,H^\prime_N)} \Big( 
\left\{ B_N  \ \leftrightarrow \ T (N, H^\prime_N) \right\}
\Big) 
\geq 1 - \exp( - C  N^{d-1} ) \, ,
\end{eqnarray}
where $H^\prime_N = 2\gd N + 2L$.

We consider also 
\begin{eqnarray*}
B^{(K)}_N = \big\{x \in B_N, \qquad x = 2Kj \quad \text{for} \;
j \in \bbZ^d \big\} \, .
\end{eqnarray*}
Then $B_N$ is covered by $\bigcup_{x \in B^{(K)}_N} b^d_K (x)$ and
\begin{eqnarray*}
\big\{ B_N \not \lra T(N,H^\prime_N) \big\}
= \bigcap_{x \in B^{(K)}_N}
\left\{ b^d_K (x)  \not \lra T(N,H^\prime_N) \right\} \, .
\end{eqnarray*}

These events are decreasing, therefore 
\eqref{eq: intermediaire 0} and FKG inequality imply
\begin{eqnarray}
\label{eq: connection 1}
\exp( - C  N^{d-1} ) \geq \prod_{x \in B^{(K)}_N}
\FKm{{\rm f}}{\boxh^* (N,H^\prime_N)} 
\left(
\left\{ b^d_K (x)  \not \lra T(N,H^\prime_N) \right\}
\right) \, .
\end{eqnarray}

Applying again the FKG property, we see that for any $x$ in
$B^{(K)}_N$
\begin{eqnarray*}
\FKm{{\rm f}}{\boxh^* (N,H^\prime_N)} 
\left(
\left\{ b^d_K (x)  \not \lra T(N,H^\prime_N) \right\}
\right) 
&\geq& 
\FKm{{\rm f}}{\boxh^* (2N,H^\prime_N)} 
\left(
\left\{ b^d_K (0)  \not \lra T(N,H^\prime_N) -x  \right\}
\right) \\
&\geq& 
\FKm{{\rm f}}{\boxh^* (2N,H^\prime_N)} 
\left(
\left\{ b^d_K (0)  \not \lra T(2N,H^\prime_N)   \right\}
\right) \, .
\end{eqnarray*}
Plugging this inequality into \eqref{eq: connection 1} leads to
\begin{eqnarray*}
- C \frac{N^{d-1}}{|B^{(K)}_N |}  
\geq 
\log \FKm{{\rm f}}{\boxh^* (2N,H^\prime_N)} 
\left(
\left\{ b^d_K (0)  \not \lra T(2N,H^\prime_N)   \right\}
\right) \, ,
\end{eqnarray*}
where $|B^{(K)}_N | = (N/2K)^{d-1}$ is the cardinal of $B^{(K)}_N$.
This completes the Proposition with $H_N = H^\prime_{N/2}$.

\qed

        \subsection{The seeds}

We  introduce now the notion of seed which will be
essential to iterate \eqref{eq: connection property}.

\begin{defi}
\label{def: top seed}
\noindent
\begin{enumerate}
\item 
We say that $b^d_K(x)$ is a {\bf seed} centered at site $x \in \bbZ^d$ if 
all the bonds lying in $b^d_K(x)$ are open.
\item 
For any $(\ell,h)$, let $\cC_K(\ell,h)$ be the event that 
$b^d_K (0)$ is connected by a path of open bonds strictly within 
$\boxh(\ell,h)$ to a seed $b^d_K(x)$ included in $T(\ell,h)$. 
\item
For any $(\ell,h)$, $\cC_K^{i}(\ell,h)$ denotes the event that $b^d_K (0)$ 
is connected strictly within $\boxh(\ell,h)$ 
to a seed $b^d_K(x)$ included in $T(\ell,h)$ and centered at
site $x$ in $T_i (\ell,h)$, with $i \leq 2^{d-1}$.
\end{enumerate}
\end{defi}

In the following, we choose $\bond \in \Theta_q$  and
$H_N$ according to Proposition \ref{pro: seed 0}.
\begin{pro}
\label{pro: seed}
For any integer $K$ there exists an integer 
$M = M(K)$ and a height 
$h = h(K,N) \in [H_N-M,  H_N]$ such that uniformly over $N$
large enough
\begin{eqnarray*}
\FKm{p, {\rm f}}{\boxh^* (N,H_N)} \Big( \cC_K (N,h) \Big) 
\geq 1 - 4 \exp \big( - C  K^{d-1} \big) \, ,
\end{eqnarray*}
for some positive constant $C$.
\end{pro}

\noindent
{\it Proof}.

For a given height $h$, define $Y(\ell,h)$ as the number of sites in 
 $T(\ell,h)$ which are connected strictly
within $\boxh (\ell,h)$ to $b^d_K (0)$ by paths of open bonds.

Let $m,M$ be two integers which are going to be chosen later.
The sequence of random heights $\cH_i$ is  defined recursively.
We set $\cH_0 = H_N -M$ and
\begin{eqnarray*}
\cH_{i+1} = \inf \big\{h > \cH_i, \quad 1 \leq Y(N,h) \leq m \big\}
\; \wedge \; H_N \, .
\end{eqnarray*}

We first check that there exists $r_\bond = r_\bond (m) <1$ such that
uniformly in $N$ and $M$
\begin{eqnarray}
\label{eq: exp decay}
\forall n \in [2,M], \qquad
\Perc_{\boxh^* (N,H_N)}^{{\rm f}} \Big(  \cH_n < H_N \Big) 
\leq r_\bond^{n-1} \, .
\end{eqnarray}

To see this, we write for $n \geq 2$
\begin{eqnarray*}
\Perc_{\boxh^*  (N,H_N)}^{{\rm f}} \Big(  \cH_n < H_N \Big) 
= \sum_{h < H_N -1}
\Perc_{\boxh^*  (N,H_N)}^{{\rm f}} \Big(  \cH_{n-1} = h, \cH_n < H_N \Big) 
 \, .
\end{eqnarray*}
Conditionally to the event $\{\cH_{n-1} = h \}$ any bond configuration
in $\{\cH_n < H_N\}$ is such that $Y(N,h+1) \geq 1$.
On the other hand, there are at most $m$ sites of $T (N,h)$ connected
to $b^d_K (0)$. Thus 
\begin{eqnarray*}
\Perc_{\boxh^*  (N,H_N)}^{{\rm f}} \Big(   \cH_n < H_N \, \Big| \, \cH_{n-1} = h\Big) 
&\leq&  
1 - \Perc_{\boxh^*  (N,H_N)}^{{\rm f}} \Big( Y(N,h+1) = 0 \, \Big| \, 
\cH_{n-1} = h\Big) 
\, , \\
&\leq&  1 - c_p^m  = r_\bond < 1 \, ,
\end{eqnarray*}
where $c_p^m >0$ is a lower bond of the energetic cost for closing the bonds
which could connect $b^d_K(0)$ to the height
$h+1$.
By iterating the procedure, we obtain \eqref{eq: exp decay}.

\vskip.3cm

For any $K$, there exists $n = n(K,m)$ such that
for any $M$ large enough and uniformly over $N$,
\begin{eqnarray}
\label{eq: decay lower}
\Perc_{\boxh^*  (N,H_N)}^{{\rm f}} 
\left( 
\sum_{h = H_N -M}^{H_N } 
1_{\{Y(N,h) > m \} } \geq M - n
\right)
\geq 1 - 2 \exp( - C K^{d-1}) \, ,
\end{eqnarray}
where $C$ is the constant of Proposition \ref{pro: seed 0}.

To see this, we write
\begin{eqnarray*}
\Perc_{\boxh^*  (N,H_N)}^{{\rm f}} 
\left( 
\sum_{h = H_N -M}^{H_N } 
1_{\{ Y(N,h) \leq m \} } \geq  n
\right)
&\leq& 
\Perc_{\boxh^*  (N,H_N)}^{{\rm f}} 
\left( Y(N,H_N) = 0 \right) \\
&& \qquad 
+\Perc_{\boxh^*  (N,H_N)}^{{\rm f}} 
\left( \cH_n < H_N \right) \, .
\end{eqnarray*}
The first term in the LHS means that there is no connection from 
$b^d_K(0)$ to the top face of $\boxh(N,H_N)$ thus it is bounded 
by Proposition \ref{pro: seed 0}. 
Let $n = n(K,m)$ be such that $r_\bond^{n-1} \leq \exp( - C K^{d-1})$. 
Then the second term is bounded by the estimate \eqref{eq: exp decay}.
Thus \eqref{eq: decay lower} is satisfied.

\vskip.3cm

We need to check that there exists $h \in [ H_N -M , H_N]$ such that
\begin{eqnarray}
\label{eq: attachment points}
\FKm{{\rm f}}{\boxh^*  (N,H_N)} \Big( Y(N,h) \geq m \Big) 
\geq 1 - 3 \exp \left( - C  K^{d-1} \right) \, .
\end{eqnarray}

Applying Tchebyshev inequality to the estimate \eqref{eq: decay lower}, 
we get
\begin{eqnarray*}
\big(1 - 2 \exp( - C K^{d-1}) \big)
(M -n)
\leq
\sum_{h = H_N -M}^{H_N} 
\Perc_{\boxh^*  (N,H_N)}^{{\rm f}} 
\big( Y(N,h) \geq m \big) \, .
\end{eqnarray*}

We choose $M = M(K,m)$ such that 
$n(K,m) < \exp(- C K^{d-1}) M$ then the above 
inequality ensures that there exists a height $h$ in $[ H_N -M ,
H_N]$ such that \eqref{eq: attachment points} holds.
The height $h$ depends on $K,N$ and $m$, nevertheless the parameter
$m$ will be fixed according to $K$ (see  \eqref{eq: borne m}).
Thus, ultimately, only the dependency on $K,N$ will remain.

\vskip.5cm

The estimate \eqref{eq: attachment points} implies that with high
probability, there are at least $m$ sites in $T(N,h)$
connected to $b^d_K(0)$. 
Thus it is enough to find at least a seed lying on top of 
one of these sites.
The probability that uniformly wrt the boundary conditions
there exists a seed of side length $K$ is
bounded from below by $c_p^{K^{d-1}}$, with $c_p>0$.
Let us choose $m = m(K)$ such that
\begin{eqnarray}
\label{eq: borne m}
\big( 1 - c_p^{K^{d-1}} \big)^{m / (4K)^{d-1}} 
< \exp( -C K^{d-1}) \, .
\end{eqnarray}

Conditioning wrt the bond configurations below $h$ which belong
to the event  $\{Y(N,h) \geq m \}$, we can find at least $m / (4K)^{d-1}$ 
sites at distance larger than $4K$ from each other.
Then the probability that a seed (lying in $T(N,h)$)
is attached on top of one of these
sites is dominated by independent Bernoulli variables. 
As $m$ has been chosen large enough, there must
be at least a seed on one of the attachment site with probability 
larger than  $1- \exp( -C K^{d-1})$.

This concludes the proof of Proposition \ref{pro: seed}.  \qed

\vskip1cm

Proposition \ref{pro: seed} has to be enhanced
in order to control the position of the seeds.

\begin{pro}
\label{pro: seed refined}
For any large $K$ there exists an integer $M = M(K)$ and a height 
$h = h(K,N) \in [H_N-M,  H_N]$ such that uniformly over $N$ large 
enough
\begin{eqnarray*}
\forall i \leq 2^{d-1}, \qquad 
\FKm{p, {\rm f}}{\boxh^* (N,H_N)} \Big( \cC_K^i (N,h) \Big) 
\geq 1 - \exp \big( - C  K^{d-1} \big) \, ,
\end{eqnarray*}
for some positive constant $C$.
The event $\cC_K^i (N,h)$ was introduced in Definition
\ref{def: top seed}.
\end{pro}

\vskip.3cm

\noindent
{\it Proof.}

Following \cite{BGN}, we use the symmetry of the model to write
\begin{eqnarray*}
&& \left(\FKm{{\rm f}}{\boxh^*  (N,H_N)} \Big( (\cC_K^i (N,h))^c \Big)
\right)^{2^{d-1}} 
=
\prod_j
\FKm{{\rm f}}{\boxh^*  (N,H_N)} \Big( (\cC_K^j (N,h))^c \Big) \\
&& \qquad \qquad  \qquad \leq
\FKm{{\rm f}}{\boxh^*  (N,H_N)} \Big( \bigcap_j (\cC_K^j (N,h))^c \Big) 
= \FKm{{\rm f}}{\boxh^*  (N,H_N)} \Big( (\cC_K (N,h))^c \Big) 
 \, ,
\end{eqnarray*}
where we used the FKG inequality and the fact that the
events $\cC_K^j (N,h)$ are increasing.

Using the lower bound derived in Proposition \ref{pro: seed},
we obtain that
\begin{eqnarray*}
\forall i \leq 2^{d-1}, \qquad
\FKm{{\rm f}}{\boxh^*  (N,H_N)} \Big( (\cC_K^i (N,h))^c \Big)
\leq  4^{1/2^{d-1}} \, \exp
\left( - \frac{C}{2^{d-1}}  K^{d-1} \right)
\, .
\end{eqnarray*}
This completes the proof. \qed

\vskip1cm

\section{The coarse graining: first renormalization step}
\label{sec:  first renormalization step}

Recall that $\bond$ is chosen in the set $\Theta_q$ 
introduced before Corollary \ref{cor: connection}.

\subsection{Top connections}

The next step is to iterate Proposition \ref{pro: seed refined}
in order to obtain a connection in a box $\boxh(N,L_N)$ for some 
$L_N \geq 3N$.

\begin{thm}
\label{thm: block type 1}
For any $\eta >0$, one can find $K = K(\eta)$  such that 
uniformly in $N$ (large enough), there exists an integer $M = M(K)$ 
and a height 
$\gp = \gp (K,N) \in [L_N - M , L_N]$ such that
\begin{eqnarray*}
\forall i \leq 2^{d-1}, \qquad
\FKm{p, {\rm f}}{\boxh^* (N,L_N)} \Big( \cC_K^i (N,\gp) \Big) 
\geq 1 - \eta \, ,
\end{eqnarray*}
where $L_N = L_N (K)$ is such that $3N \leq L_N - M(K)
\leq \gp (K,N)$ (see \eqref{eq: parametre L_N}).
\end{thm}

\vskip.3cm

\noindent
{\it Proof.}

Let $H_N = \gd N + L$, where $\gd$ and $L$ are chosen according to
Proposition \ref{pro: seed 0}.
Fix $K$ large enough such that
\begin{eqnarray*}
\left( 1 - \exp \big( - C  K^{d-1} \big) 
\right)^{1 + 6 /\gd} \geq 1 - \eta \, ,
\end{eqnarray*}
where $C$ is the constant obtained in Proposition 
\ref{pro: seed refined}.
Then choose $N$ large enough and set $h = h(K,N), M= M(K)$ according 
to Proposition \ref{pro: seed refined}.\\

We set 
\begin{eqnarray}
\label{eq: parametre L_N}
n = 6/\gd, \qquad L^\prime_N = H_N + n h(K,N), \qquad
L_N = L^\prime_{N/2} \, .
\end{eqnarray}
Let $\cA_i$ be the event that $b^d_K(0)$ is connected to a 
seed in $T \big( N,i h \big)$ by an open path strictly 
contained within $\boxh (2N , i \, h)$
\begin{eqnarray*}
\forall i \leq n, \qquad
\cA_i = 
\left\{
b^d_K (0) \lra \seed \subset T \big( N,i \, h \big)
\right\}
\, .
\end{eqnarray*}
We stress the fact that the event $\cA_i$ differs from $\cC_K (2N,ih)$
since the attachment site of a seed must be at distance less than $N$ 
from the axis $\{x_1 = 0, \dots, x_{d-1} =0 \}$, nevertheless the 
open path leading to this seed is allowed to use any bonds within 
$\boxh (2N , i \, h)$.
If $\gp^\prime(K,N) = (n+1) h(K,N)$ then
$\cA_{n+1} \subset \cC_K (2N, \gp^\prime)$.

\begin{figure}
\begin{center}
\leavevmode
\epsfysize = 7 cm
\psfrag{0}{$b^d_K(0)$}
\psfrag{1}[Bc]{$y^{(1)}$}
\psfrag{2}{$y^{(2)}$}
\psfrag{3}{$y^{(3)}$}
\psfrag{4}{$y^{(4)}$}
\psfrag{5}{$y^{(5)}$}
\psfrag{6}[Br]{{\small $y^{(6)}$}}
\psfrag{n}{$N$}
\psfrag{m}{$2N$}
\psfrag{h}[Br]{$h(K,N)$}
\psfrag{L}{$L^\prime_N$}
\epsfbox{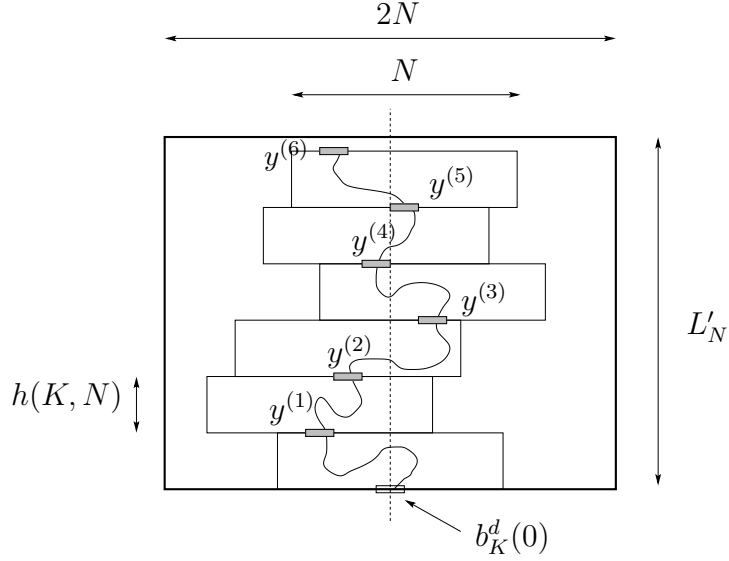}
\end{center}
\caption{The seeds are connected according to the steering rule.}
\end{figure}

Conditioning wrt $\pi_n$ the bond configuration below the
 height $n h(K,N)$, we get
\begin{eqnarray*}
&& \FKm{{\rm f}}{\boxh^* (2N, L^\prime_N)} \left( \cA_{n+1} \right)
\geq 
\FKm{{\rm f}}{\boxh^* (2N, L^\prime_N)} 
\left( \cA_n \cap \cA_{n+1}  \right)\\
&& \qquad \geq
\FKm{{\rm f}}{\boxh^* (2N, L^\prime_N)} 
\left( \cA_n \; 
\FKm{\pi_n, {\rm f}}{\boxh^* (2N, L^\prime_N)} ( \cA_{n+1})
\right)\\
&& \qquad \geq 
\FKm{{\rm f}}{\boxh^* (2N, L^\prime_N)} 
\left( \cA_n  \; 
\FKm{{\pi_n, \rm f}}{\boxh^* (2N,H_N)} 
\left( \big\{ b_K^d(y^{(n)}) \lra \seed \subset T(N, (n+1)h) \big\} 
\right) \right) \, ,
\end{eqnarray*}
where $y^{(n)}$ is defined for each $\pi_n$ as the attachment site 
of a seed in $T(N, n h)$ 
(if there are many such sites choose the smallest one 
wrt the lexicographic order).
The box $y^{(n)} + \boxh(N,H_N)$ is included in $\boxh(2N,H_N)$ and
with high probability $b^d_K (y^{(n)})$ is connected to a seed in
$y^{(n)} + T(N, h)$. Nevertheless,
piling up the boxes in an arbitrary way would give little control 
on the location of the boxes and therefore the connections might not
remain in the box $\boxh(2 N, L^\prime_N)$. 
In order to ensure that this seed is in $T(N,(n+1)h)$ one needs to
choose a proper subfacet according to the {\bf steering rule} described below.

We write $y^{(n)} = \{ y^{(n)}_i \}_{i \leq d}$.
Choose $j(\pi_n)$ such that the site $\{-y^{(n)}_1, \dots, 
-y^{(n)}_{d-1},h\}$ belongs to $T_{j(\pi_n)} (N,h)$.
By construction $y^{(n)} + T_{j(\pi_n)} (N,h) \subset T(N,(n+1) h)$.
\begin{eqnarray*}
\FKm{{\pi_n, \rm f}}{\boxh^* (2N,H_N)} 
\left( \big\{ b_K^d(y^{(n)}) \lra \seed \subset T(N, (n+1)h) \big\} \right)
\geq 
\FKm{{\pi_n, \rm f}}{\boxh^* (2N,H_N)} 
\left( \cC^{j(\pi_n)}_K (N,h) \right)
\, ,
\end{eqnarray*}
where the event $ \cC^{j(\pi_n)}_K (N, h)$ is a translate of the event
introduced in Definition \ref{def: top seed} and is supported by
the domain $y^{(n)} + \boxh^* (N,H_N)$.
As $\cC^{j}_K (N, h)$ is an increasing event, we derive
a lower bound uniformly over the configurations $\pi_n$ in $\cA_n$
\begin{eqnarray*}
\FKm{\pi_n, {\rm f}}{\boxh^* (2 N,H_N)} 
\left( \cC^{j(\pi_n)}_K (2N,h)  \right)
\geq 
\FKm{{\rm f}}{\boxh^* (N,H_N)} 
\left( \cC^1_K (N,h)\right)
\geq 1 - \exp \big( -C K^{d-1} \big) \, ,
\end{eqnarray*}
where the last estimate follows from Proposition 
\ref{pro: seed refined}.

Proceeding recursively, and applying the steering rule at each
step we obtain 
\begin{eqnarray*}
\FKm{{\rm f}}{\boxh^* (N,L_N)} \Big( \cC_K (N,\gp) \Big) 
\geq 1 - \eta \, ,
\end{eqnarray*}
with $L_N = L^\prime_{N/2}$ and $\gp(K,N)= \gp^\prime(K,N/2)$. 
Using the same strategy as in Proposition \ref{pro: seed refined},
the position of the seeds can be controlled. 
This completes the proof.

\qed

\medskip

The steering rule was  introduced in \cite{BGN}.
It will be also  extremely useful to perform the dynamic 
renormalization in a two-dimensional slab (see Section
\ref{sec: second renormalization level}).

\vskip.5cm

\begin{rem}
\label{rem: m}
The previous Theorem enables us to generalize  inequality
\eqref{eq: attachment points} to the domain $\boxh(N,L_N)$.
For a given $\eta>0$, choose $K$ and $\gp(K,N)$ according to
Theorem \ref{thm: block type 1}.
Then, for any $m>0$, there exists $M^\prime = M^\prime (\eta,m)$
such that for any $N$ large enough, one can find a height
$\hat \gp = \hat \gp (\eta,m,K,N) \in [\gp(K,N) - M^\prime ,\gp(K,N)
]$ so that
\begin{eqnarray}
\label{eq: random height brick}
\FKm{{\rm f}}{\boxh^* (N,L_N)} \big( Y(N, \hat \gp ) > m \big) 
\geq 
1 - 2 \eta \, .
\end{eqnarray}
The derivation of this inequality follows the scheme of the 
proof of \eqref{eq: attachment points}.
\end{rem}

\vskip.5cm

        \subsection{Occupied blocks}

Theorem \ref{thm: block type 1} implies that a connection occurs 
with high probability from a small region around the origin to
a height almost at the top of a box $\boxh(N,L_N)$.
To perform the dynamic renormalization, connections
from the bottom to the lateral sides of a box will be necessary.
This will require different arguments.\\

We state now the counterpart of Definition
\ref{def: top seed} for the lateral connections.

\begin{defi}
\label{def: lateral seed}
\noindent
\begin{enumerate}
\item
For any integers $(\ell, h)$,
denote by $\hat \cC_K(\ell, h)$ the event that $b_K^d (0)$ is connected 
by open bonds strictly within $\boxh(\ell,h)$ to a seed lying entirely in
$S(\ell,h)$ (the union of the sides).
Depending on the orientation of the side, the corresponding
seed will be of the type $b^j_K (\cdot)$ for some appropriate $j$.
\item
We  also define the collection of events
$\hat \cC_K^i(\ell,h)$  for all $i \leq (d-1) 2^{d-1}$. 
Any bond configuration in $\hat \cC_K^{i}(\ell,h)$ is such that $b^d_K
(0)$ is connected strictly within $\bbB(\ell,h)$ to a seed 
included in $S(\ell,h)$ and centered at a site in $S_i (\ell,h)$. 
\item
A block $\bbB(\ell,h)$ is said {\bf occupied} if the
bond configuration within this block belongs to the intersection
of all the sets $\cC_K^{i}(\ell,h)$ and $\hat \cC_K^{j}(\ell,h)$. 
\end{enumerate}
\end{defi}

\vskip.5cm

The main step will be to prove that there is a class of blocks 
which are {\bf occupied} with high probability.

\begin{thm}
\label{thm: block type 2}
Fix $\eta>0$, then there exists $K = K(\eta)$ such that for $N$ 
large enough one can find $M(K)$,
$h = h(K,N) \in [L^\prime_N- M(K), L^\prime_N]$ and 
$\ell = \ell(K,N) \leq h (K,N)$ so that
\begin{eqnarray*}
\FKm{p, {\rm f}}{\boxh^* (N,L^\prime_N)} 
\left( \bigcap_i  \cC_K^{i}(\ell,h) \ \bigcap_j \hat \cC_K^{j}(\ell,h)  \right) 
\geq 1 - \eta \, ,
\end{eqnarray*}
with $L^\prime_N = L^\prime_N(K) \geq 3N$.
\end{thm}

\vskip.5cm

The Theorem is a consequence of Propositions 
\ref{pro: connections laterales} and \ref{pro: top connections}
which are derived in the following Subsections.

        \subsubsection{Lateral connections}

Recall that $Y(\ell, h)$ denotes the number of sites in 
$T(\ell, h)$  which are connected strictly
within $\boxh (\ell, h)$ to  $b^d_K (0)$ by open bonds.
We also define $X(\ell, h)$ as the number of sites in the sides
$S(\ell, h)$  which are connected strictly
within $\boxh (\ell, h)$ to  $b^d_K (0)$ by a path of open bonds.

\vskip.5cm

\begin{pro}
\label{pro: connections laterales}
Fix $\eta>0$, then there exists $K = K(\eta)$ such that 
for $N$ large enough one can find 
$\gp^1 = \gp^1(K,N)$, $\ell^1 = \ell^1(K,N)$  and 
$L^\prime_N \geq 3N$ so that
\begin{eqnarray*}
\forall j \leq (d-1)2^{d-1}, \qquad
\FKm{{\rm f}}{\boxh^* (N,L^\prime_N)} 
\big(  {\hat \cC}^j_K (\ell^1,\gp^1)   \big) 
\geq 1 - \gep_1 (\eta) \, ,
\end{eqnarray*}
where $\gep_1 (\cdot)$ converges to 0 as $\eta$ tends to 0. 
\end{pro}

\noindent
{\it Proof.}

We start by fixing the parameters according to $\eta$.
Remark \ref{rem: m} implies that there exists $K = K(\eta)$
and a sequence $L_n$, such that for any $m>0$ one can find
$M = M(\eta,m,K)$ and $\gp_n = \gp_n(\eta,m,K) \in [L_n -M, L_n]$ 
such that for any $n$ large enough,
\begin{eqnarray}
\label{eq: parametres}
\FKm{{\rm f}}{\boxh^* (n,L_n)} \big( Y(n, \gp_n ) \geq  m \big) 
\geq 1 -  \frac{1}{2} (\eta^2)^{2^{d-1}}  \gg 1 -  \eta^2 \, .
\end{eqnarray}
Finally, let $m = m(\eta,K)$ be such that (see \eqref{eq: borne m})
\begin{eqnarray*}
\big( 1 - c_p^{K^{d-1}}  \big)^{m / (4K)^{d-1}}
< \frac{1}{2} (\eta^2)^{2^{d-1}} \, .
\end{eqnarray*}

As $m$ has been chosen large enough, \eqref{eq: parametres}
implies that uniformly over $n$ large enough (see Propositions
\ref{pro: seed} and \ref{pro: seed refined})
\begin{eqnarray}
\label{eq: parametres 2}
\forall i \leq 2^{d-1}, \qquad 
\FKm{{\rm f}}{\boxh^* (n,L_n)} \big( \cC^i_K (n, \gp_n ) \big) 
\geq 1 -   \eta^2 \, .
\end{eqnarray}

\medskip

For any $n$ large enough, there exists $\ell_n$ such that 
\begin{eqnarray}
\label{eq: height step 1}
\begin{cases}
\FKm{{\rm f}}{\boxh^* (n,L_n)} 
\big(  Y(\ell_n - 1,\gp_n) \leq m -1 \big) > \eta \, , \\
\FKm{{\rm f}}{\boxh^* (n,L_n)} 
\big(  Y(\ell_n,\gp_n) \geq m  \big) \geq 1- \eta \, .
\end{cases}
\end{eqnarray}
To see that \eqref{eq: height step 1} is well defined, first notice that
\begin{eqnarray*}
\FKm{{\rm f}}{\boxh^* (n,L_n)} 
\Big(  Y(\ell -1,\gp_n) \leq  m -1  \Big) = 1 - 
\FKm{{\rm f}}{\boxh^* (n,L_n)} 
\Big(  Y(\ell -1,\gp_n) \geq  m \Big) \, .
\end{eqnarray*}
Then it is enough to
observe that the function $\ell \to \FKm{{\rm f}}{\boxh^* (n,L_n)} 
\big(  Y(\ell,\gp_n) \geq m  \big)$ is non-decreasing and
\begin{eqnarray*}
\FKm{{\rm f}}{\boxh^* (n,L_n)} 
\Big(  Y(1,\gp_n) \geq m  \Big) = 0
\quad \text{and} \quad
\FKm{{\rm f}}{\boxh^* (n,L_n)} 
\Big(  Y(n,\gp_n) \geq m  \Big) \geq 1- \eta \, .
\end{eqnarray*}
Thus one defines  $\ell_n$ as the smallest value
for which \eqref{eq: height step 1} holds.\\

The proof will be decomposed into 3 steps.\\

\noindent
{\bf Step 1.}

In this first step, we show that  lateral connections occur in a domain
smaller than $\bbB(\ell_n, \gp_n)$.  
\begin{lem}
\label{lem: lateral seed}
There exists an integer $U = U(\eta)$ such that for any $n$ large 
enough, one can find $\ell^\prime_n \in [\ell_n -1 -U, \ell_n - 1]$ and
\begin{eqnarray*}
\FKm{{\rm f}}{\boxh^* (n,L_n)} \Big(\hat \cC_K (\ell^\prime_n ,\gp_n) \Big) 
\geq 1 - 2 \eta \, ,
\end{eqnarray*}
where $K$ and $\gp_n$ were chosen in \eqref{eq: parametres} and $\ell_n$ was
determined in \eqref{eq: height step 1}.
\end{lem}

\noindent
{\it Proof}.

We first remark  that for $\ell_n$ defined as above then 
\begin{eqnarray}
\label{eq: Z control}
\FKm{{\rm f}}{\boxh^* (n,L_n)} 
\Big(  X(\ell_n-1,\gp_n) = 0  \Big) \leq \eta \, .
\end{eqnarray}
To see this, we apply \eqref{eq: parametres} and obtain
\begin{eqnarray*}
&& \eta^2 \geq
\FKm{{\rm f}}{\boxh^* (n,L_n)} \big( Y( n,\gp_n ) <  m \big) \\
&& \qquad \qquad 
\geq 
\FKm{{\rm f}}{\boxh^* (n,L_n)} \big( Y( \ell_n-1, \gp_n ) \leq m -1, \;
X(\ell_n-1,\gp_n) = 0 \big) \\
&& \qquad \qquad 
\geq 
\FKm{{\rm f}}{\boxh^* (n,L_n)} \big( Y( \ell_n-1, \gp_n ) \leq m -1\big) 
\;
\FKm{{\rm f}}{\boxh^* (n,L_n)} \big( X(\ell_n-1,\gp_n) = 0 \big) 
\, ,
\end{eqnarray*}
where the last inequality follows from the fact that both events
in the RHS are decreasing. 
Thus inequality \eqref{eq: height step 1} implies
\eqref{eq: Z control}.\\

Let us also check that $\ell_n$ diverges with $n$. 
First notice that there exists $c_p > 0$ such that
\begin{eqnarray}
\label{eq: cut step 1}
\forall k \leq \gp_n, \qquad
\FKm{{\rm f}}{\boxh^* (n,L_n)} \big( Y( \ell_n ,k ) = 0 \; \big| 
\; Y( \ell_n , k-1 ) \geq 1 \big) 
\geq c_p^{\ell_n} \, ,
\end{eqnarray}
where $c_p^{\ell_n}$ is simply the cost of closing all the bonds 
linking the level $k-1$ to the level $k$.
\begin{eqnarray*}
1 - \eta \leq
\FKm{{\rm f}}{\boxh^* (n,L_n)} \big( Y( \ell_n, \gp_n ) \geq m  \big) 
\leq 
\FKm{{\rm f}}{\boxh^* (n,L_n)} \big( \forall k \leq \gp_n, \quad
 Y( \ell_n , k ) \geq 1 \big)  \, .
\end{eqnarray*}
In order to estimate the RHS, we follow the strategy used in 
Proposition \ref{pro: seed}.
We proceed recursively by conditioning at each level and using 
\eqref{eq: cut step 1}.
\begin{eqnarray*}
1 - \eta \leq (1 - c_p^{\ell_n} )^{\gp_n} \leq \exp( - c_p^{\ell_n} \,
n)  \, ,
\end{eqnarray*}
where we used that $\gp_n \geq 3 n$.
This implies that for $n$ large enough 
\begin{eqnarray}
\label{eq: ell_n}
\ell_n \geq \frac{1}{\log ( 1/c_p)} \log n \, .
\end{eqnarray}

\medskip

The derivation of Lemma \ref{lem: lateral seed} 
follows closely the one of Proposition \ref{pro: seed}.
Let $u,U$ be two integers which are going to be fixed later.
Define recursively the sequence of random lengths by
$\cL_0 = \ell_n - 1 - U$ and
\begin{eqnarray*}
\cL_{i+1} = \inf \big\{ l > \cL_i, \quad 1 \leq X(l,\gp_n) \leq u \big\}
\; \wedge \; (\ell_n - 1) \, .
\end{eqnarray*}
Notice that for $n$ large enough $\cL_0$ is positive (see 
\eqref{eq: ell_n}).

As in \eqref{eq: exp decay}, we get
\begin{eqnarray*}
\forall i \leq U, \qquad
\FKm{{\rm f}}{\boxh^* (n,L_n)} \Big(  \cL_i < \ell_n -1 \Big) 
\leq r_\bond (u)^{i-1} \, ,
\end{eqnarray*}
where $r_\bond(u) < 1$ is related to the probabilistic cost of closing
at most $u$ open bonds (see \eqref{eq: exp decay}).

Fix $u = u(\eta)$ large enough such that 
$\big( 1 - c_p^{K^{d-1}} \big)^{u / (4K)^{d-1}} < \eta/3$
(see \eqref{eq: borne m}).
Then choose $v = v(\eta,u)$ such that $r_\bond (u)^{v-1} < \eta/3$.
Following the derivation of \eqref{eq: decay lower},
we get
\begin{eqnarray*}
\FKm{{\rm f}}{\boxh^*  (n,L_n)}
\left( 
\sum_{l = \ell_n - 1 - U}^{\ell_n - 1} 
1_{\{ X(l,\gp_n) \leq u \} } \geq  v
\right)
&\leq& 
\FKm{{\rm f}}{\boxh^*  (n,L_n)}
\left( X(\ell_n - 1,\gp_n)  = 0 \right) \\
&& \qquad 
+ \FKm{{\rm f}}{\boxh^*  (n,L_n)}
\left( \cL_v < \ell_n - 1  \right) 
\leq \frac{4}{3} \eta \, ,
\end{eqnarray*}
where we used \eqref{eq: Z control}.

Finally as in the derivation of \eqref{eq: attachment points},
one can find $U(\eta,u)$ large enough such that uniformly in $n$
there exists $\ell^\prime_n \in [ \ell_n - 1 -U  , \ell_n - 1]$ 
\begin{eqnarray*}
\FKm{{\rm f}}{\boxh^* (n,L_n)} \Big( X(\ell^\prime_n,\gp_n) \geq u \Big) 
\geq 1-  \frac{5}{3} \eta \, .
\end{eqnarray*}

Therefore the probability that a seed lies on top of one of the
$u$ attachment sites is large.
This concludes  Lemma \ref{lem: lateral seed}.

\qed

\vskip.5cm

\begin{lem}
\label{lem: lateral seed improved}
There exists an integer $U = U(\eta)$ such that for any $n$ large
enough, one can find $\ell^\prime_n \in [\ell_n -1 -U, \ell_n - 1]$ and
\begin{eqnarray}
\label{eq: connection bords}
\forall i \leq (d-1) 2^{d-1}, \qquad
\FKm{{\rm f}}{\boxh^* (n,L_n)} \Big(\hat \cC^i_K (\ell^\prime_n ,\gp_n) \Big) 
\geq 1 - (2 \eta)^{\frac{1}{(d-1) 2^{d-1}}} \, ,
\end{eqnarray}
where $K$ and $\gp_n$ were chosen in \eqref{eq: parametres} and $\ell_n$ was
determined in \eqref{eq: height step 1}.
\end{lem}

The derivation is similar to the one of Proposition \ref{pro: seed refined}.

\vskip.8cm

\noindent
{\bf Step 2.}
At this stage, we know that with high probability $b^d_K(0)$ is connected
to a seed in each side of the box $\boxh(\ell^\prime_n,\gp_n)$ with
$\ell^\prime_n \in [\ell_n -1 -U, \ell_n -1]$.
As $\boxh(\ell^\prime_n, \gp_n)$ is contained in $\boxh(\ell_n, \gp_n)$, 
we do not
have any information on the connections from $b^d_K(0)$ to the top face
$T(\ell^\prime_n,\gp_n)$, thus we need to prove that the lateral connections 
occur as well in a larger box.\\

\begin{figure}
\begin{center}
\leavevmode
\epsfysize = 5 cm
\psfrag{p1}[c]{$\gp_n$}
\psfrag{p2}[c]{$\gp_{t_n}$}
\psfrag{LN}[r]{$L^\prime_N$}
\psfrag{e1}{$\vec{e}_1$}
\psfrag{ed}{$\vec{e}_d$}
\psfrag{0}{$b^d_K (0)$}
\psfrag{z1}{$z^{(1)}$}
\psfrag{z2}{$z^{(2)}$}
\psfrag{z3}{$z^{(3)}$}
\psfrag{B}{$\bbB (\ell_n, \gp_n)$}
\epsfbox{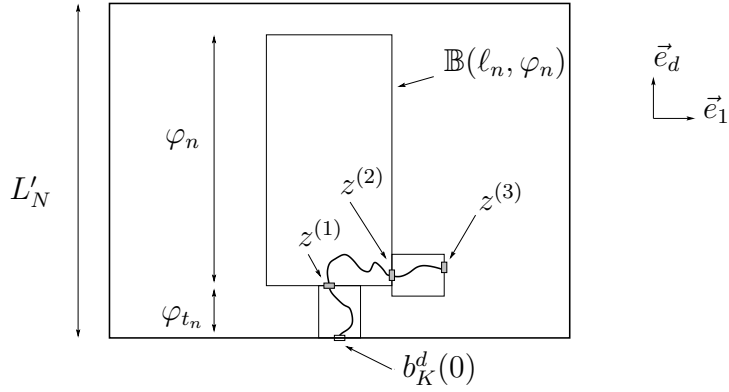}
\end{center}
\caption{The events $\cE^1,\cE^2$ and $\cE^3$ are depicted.}
\label{fig: occupied}
\end{figure}

We will use $n$ as the reference scale and consider blocks
defined in terms of the new parameters (see figure \ref{fig: occupied})
\begin{eqnarray}
\label{eq: variables 1}
&& t_n = [\sqrt{ \log  n}], \quad
\gp^1_n = t_n + \gp_n +\gp_{t_n} , \\ 
&& N =  n + t_n +  L_{t_n}, \quad
L^\prime_N = \gp_n + 4 \gp_{t_n} + L_{t_n}  \, . \nonumber
\end{eqnarray}
As $\gp_n \geq 3n$ (see Theorem \ref{thm: block type 1}), one has
$L^\prime_N \geq 3N$.

We are going to show that $b^d_K(0)$ is connected with high
probability strictly within $\bbB(\ell^\prime_n + 3 t_n,\gp^1_n)$
to the sides of this block.
\begin{lem}
\label{lem: premiere connection}
Uniformly over $n$ large enough,
\begin{eqnarray}
\label{eq: lateral 0}
\FKm{{\rm f}}{\boxh^*  (N,L^\prime_N)}
\big(  X( \ell^\prime_n + 3 t_n, \gp^1_n ) \geq 1 \big)
\geq 1 - \gd( \eta) \, ,
\end{eqnarray}
where  $\lim_{\eta \to 0} \gd(\eta) = 0$.
\end{lem}

        \begin{proof}

Let $\cE^1 = \cC^1_K(t_n, \gp_{t_n})$ be the event that $b^d_K (0)$ 
is connected to a seed attached in $T^1( t_n, \gp_{t_n})$ by a path 
strictly within $\bbB(t_n, \gp_{t_n})$. 
From \eqref{eq: parametres 2} and the FKG inequality, we have
\begin{eqnarray}
\label{eq: lateral 1}
\FKm{{\rm f}}{\boxh^*  (N,L^\prime_N)}
\left(  \cE^1 \right)
\geq 1 - \eta \, .
\end{eqnarray}

For any bond configuration in $\cE^1$, we denote by $z^{(1)}$ 
the center of the seed in $T^1( t_n, \gp_{t_n})$ with the
largest coordinate in the lexicographic order.
By construction
\begin{eqnarray*}
\forall i \leq d-1, \quad 0 \leq z^{(1)}_i \leq t_n, 
\quad \text{and} \quad
z^{(1)}_d = \gp_{t_n} \, .
\end{eqnarray*}
The set $\cE^1$ can be partitioned into 
\begin{eqnarray*}
\cE^1 = \bigcup_{z^{(1)} \in T_1( t_n, \gp_{t_n})}
\cE^1 ( z^{(1)}) \, ,
\end{eqnarray*}
where $\cE^1(z^{(1)})$ is the set of bond configurations in 
$\cE^1$ for which the seed connected to $b^d_K(0)$ 
with the largest coordinate is centered in $z^{(1)}$.

\medskip

Given $z^{(1)}$ in $T^1( t_n, \gp_{t_n})$, we introduce the event
$\cE^2 (z^{(1)}) $ such that $b_K^d (z^{(1)})$ is connected strictly 
within $z^{(1)} + \bbB(\ell^\prime_n, \gp_n)$ to a seed $b_K^1 (z^{(2)})$ 
lying in $z^{(1)} + S_1 (\ell^\prime_n, \gp_n)$ 
with attachment site $z^{(2)} = (z^{(2)}_i)_{i \leq d}$,
\begin{eqnarray*}
z^{(2)}_1 = \ell^\prime_n + z^{(1)}_1, \quad 
\gp_{t_n} \leq z^{(2)}_d \leq \gp_{t_n} + \gp_n, \qquad
\forall i \in \{2,d-1\}, \quad 0 \leq z^{(2)}_i \leq t_n + \ell^\prime_n \, .
\end{eqnarray*}
Once again, $z^{(2)}$ is chosen wrt some arbitrary order and therefore
is uniquely determined.

Let $\cF^{(1)}$ be the $\gs$-algebra generated by the bond
configurations in $\bbB( t_n, \gp_{t_n})$.
Thanks to Lemma \ref{lem: lateral seed improved}, we have 
for any $z^{(1)}$,
\begin{eqnarray}
\label{eq: lateral 2}
\FKm{{\rm f}}{\boxh^*  (N,L^\prime_N)}
\left(  \cE^2 (z^{(1)}) \, \big| \, \cF^{(1)} \right)
\geq \FKm{{\rm f}}{\boxh^*  (n,L_n)}
\left(  \hat \cC^1_K ( \ell^\prime_n , \gp_n)  \right)
\geq 1 -  (2 \eta)^{\frac{1}{(d-1) 2^{d-1}}} \, ,
\end{eqnarray}
where we used the FKG property and the fact that for any
$z^{(1)} \in T^1( t_n, \gp_{t_n})$ then $z^{(1)} + \bbB(n,L_n)$
remains in $\bbB(N,L^\prime_N)$.

\medskip

Similarly, $\cE^2(z^{(1)}) $ can be partitioned into 
\begin{eqnarray*}
\cE^2(z^{(1)}) = 
\bigcup_{z^{(2)} \in \{ z^{(1)} + S_1( \ell^\prime_n, \gp_{n}) \}}
\cE^2(z^{(1)},z^{(2)}) \, .
\end{eqnarray*}
For a given $z^{(2)}$, let $\cF^{z^{(2)}}$ be the $\gs$-algebra 
generated by the bond configurations in the half plane 
$\{x, \quad x_1 \leq z^{(2)}_1 \}$.
Define the event $\cE_3(z^{(2)})$ that $b^1_K (z^{(2)})$ is connected
to a site $z^{(3)}$ which satisfies
\begin{eqnarray*}
\begin{cases}
z^{(3)}_1 = z^{(2)}_1 +  \gp_{t_n} \geq \ell^\prime_n + 3 t_n \, , \\
z^{(2)}_i -  t_n \leq z^{(3)}_i \leq z^{(2)}_i + t_n, 
\qquad \forall i \in \{2,d-1\} \, ,\\
0 \leq z^{(3)}_d \leq \gp_{t_n} + \gp_n + t_n\, ,
\end{cases}
\end{eqnarray*}
furthermore the connection must occur strictly within the set
$\{z^{(2)}_1, z^{(2)}_1 +  \gp_{t_n} \} \times
\{- t_n, \ell^\prime_n + 2 t_n \}^{d-2} \times
\{0, \gp_n + \gp_{t_n} + t_n\}$.
In the previous computation we used the fact that $\gp_n \geq 3n$
(see Theorem \ref{thm: block type 1}).\\

We are going to check that \eqref{eq: lateral 0} is implied by
the following inequality 
\begin{eqnarray}
\label{eq: lateral 3}
\FKm{{\rm f}}{\boxh^*  (N,L^\prime_N)}
\left( \cE^3 (z^{(2)}) \, \big| \, \cF^{z^{(2)}} \right)
\geq  1 - \eta \, .
\end{eqnarray}
By construction
\begin{eqnarray*}
&&\FKm{{\rm f}}{\boxh^* (N,L^\prime_N)} 
\big( X( \ell^\prime_n + 3 t_n, \gp^1_n) \geq 1  \big) \\
&&\quad \geq \sum_{z^{(1)} \in T_1( t_n, \gp_{t_n})} 
\sum_{z^{(2)} \in \{ z^{(1)} + S_1( \ell^\prime_n, \gp_{n}) \}}
\FKm{{\rm f}}{\boxh^*  (N,L^\prime_N)}
\left( \cE^1 (z^{(1)}) \cap \cE^2 (z^{(1)},z^{(2)}) \cap 
\cE^3 (z^{(2)}) \right)
\, ,
\end{eqnarray*}
where the events in the RHS are disjoint by construction.
Applying \eqref{eq: lateral 3} we obtain
\begin{eqnarray*}
\FKm{{\rm f}}{\boxh^* (N,L^\prime_N)} 
\big( X( \ell^\prime_n + 3 t_n, \gp^1_n) \geq 1 \big) 
\geq \sum_{z^{(1)} \in T_1( t_n, \gp_{t_n})} 
\FKm{{\rm f}}{\boxh^*  (N,L^\prime_N)}
\left( \cE^1 (z^{(1)}) \cap \cE^2 (z^{(1)}) \right)
(1 - \eta) \, .
\end{eqnarray*}
Using \eqref{eq: lateral 2} and then 
\eqref{eq: lateral 1}, Lemma \ref{lem: premiere connection} 
is satisfied with
\begin{eqnarray*}
\gd(\eta) = 1- \left( 1- (2 \eta)^{\frac{1}{(d-1) 2^{d-1}}} 
\right) ( 1 - \eta)^2
\leq (2 \eta)^{\frac{1}{(d-1) 2^{d-1}}} + 2 \eta \, .
\end{eqnarray*}

\vskip.5cm

It remains to derive  \eqref{eq: lateral 3}  by piling
 a block oriented in the direction
$\vec{e}_1$.
We introduce $\bbB^{(1)} (L,H) = \{0,H \}\times \{-L,L\}^{d-1}$
and its top face is $T^{(1)} (L,H) = \{H\}\times \{-L,L\}^{d-1}$.
According to \eqref{eq: ell_n},  $\ell^\prime_n \gg t_n$ for $n$
large enough.
Thus $z^{(2)} + \bbB^{(1)} (t_n,L_{t_n})$ is 
included in $\bbB(N,L^\prime_N)$ but does not intersect $\bbB(t_n, \gp_{t_n})$.
Inequality \eqref{eq: parametres 2} implies that with
probability larger than $1 - \eta$, the seed $b^1_K(z^{(2)})$ 
is connected within $z^{(2)} + \bbB^{(1)} (t_n,\gp_{t_n})$ to a seed  
in $z^{(2)} + T^{(1)} (t_n,\gp_{t_n})$ attached at the site $z^{(3)}$.




\end{proof}

\vskip.5cm

\noindent
{\bf Step 3.}

Lemma \ref{lem: premiere connection} implies that with high
probability, $b^d_K(0)$ is connected to 
$S( \ell^\prime_n + 3 t_n, \gp^1_n )$ within
$\bbB( \ell^\prime_n + 3 t_n, \gp^1_n )$.
Thus we can proceed as in the first step and derive on a 
larger domain a result similar to \eqref{eq: connection bords}.
There exists an integer $U = U(\eta)$ such that for any $n$ large 
enough, one can find $\ell^1_n$ in 
$[\ell^\prime_n + 3 t_n -U, \ell^\prime_n + 3 t_n]$ such that
\begin{eqnarray}
\label{eq: connection bords 2}
\forall j \leq (d-1) 2^{d-1}, \qquad
\FKm{{\rm f}}{\boxh^* (N,L^\prime_N)} \Big(\hat \cC^j_K (\ell^1_n ,\gp^1_n) \Big) 
\geq 1 - \big( 2 \gd(\eta) \big)^{\frac{1}{(d-1) 2^{d-1}}} \, .
\end{eqnarray}
This completes the proof of Proposition \ref{pro: connections
laterales} with $\gep_1 (\eta)$ defined by the RHS of the 
previous inequality.


\qed

	\subsubsection{Top connections}

For any $\eta$, the parameters are fixed according to Proposition 
\ref{pro: connections laterales}. In particular, $K$ was introduced in
\eqref{eq: parametres},  $\gp^1_n, N$ in \eqref{eq: variables 1} and 
$\ell^1_n$ in \eqref{eq: connection bords 2}.

\begin{pro}
\label{pro: top connections}
There exists $\ell^2 = \ell^2(K,N)$ and $\gp^2 = \gp^2 (K,N)
 \in [L^\prime_N - M(K), L^\prime_N]$  such that
\begin{eqnarray*}
\forall i \leq 2^{d-1}, \qquad
\FKm{{\rm f}}{\boxh^* (N,L^\prime_N)} 
\big(  \cC^i_K (\ell^2,\gp^2)   \big) 
\geq 1 - \gep_2(\eta) \, ,
\end{eqnarray*}
where $\gep_2 (\cdot)$ converges to 0 as $\eta$ tends to 0. 
By construction, $\ell^2(K,N)\leq \ell^1 (K,N)$ and 
$\gp^2 (K,N)\geq \gp^1 (K,N)$.
\end{pro}

\begin{proof}
Let $\ell_n$ be defined according to \eqref{eq: height step 1}.
Proceeding as in the derivation of \eqref{eq: parametres} and 
since  $m$ has been chosen large enough, we have
\begin{eqnarray}
\label{eq: houba}
\FKm{{\rm f}}{\boxh^* (n,L_n)} 
\big(  \cC_K(\ell_n,\gp_n)  \big) \geq 1- 2 \eta \, .
\end{eqnarray}
We recall that \eqref{eq: parametres} also implies for $n$ 
large enough that
\begin{eqnarray}
\label{eq: houba2}
\forall i \leq 2^{d-1}, \qquad
\FKm{{\rm f}}{\boxh^* (t_n,\gp_{t_n})} 
\big(  \cC^i_K (t_n,\gp_{t_n})  \big) \geq 1-  
\eta^2  \, .
\end{eqnarray}

As in \eqref{eq: variables 1}, the scaling parameter is $n$
\begin{eqnarray}
\label{eq: variables 2}
\ell^2(K,N) = \ell^2_n = \ell_n + t_n, \qquad 
\gp^2 (K,N) = \gp^2_n = 5 \gp_{t_n} + \gp_n \, .
\end{eqnarray}
According to Theorem \ref{thm: block type 1},
$\gp_{t_n} \geq 3 t_n$, thus for $n$ large enough 
\begin{eqnarray*}
\ell^1_n = \ell_n -1 - 2U + 3 t_n \geq \ell^2_n , \qquad 
\gp^1_n = t_n  + \gp_n + \gp_{t_n} \leq  \gp^2_n \, .
\end{eqnarray*}
Finally, we recall that $L^\prime_N =\gp_n + 4 \gp_{t_n} + L_{t_n}$ 
(see \eqref{eq: variables 1})
so that $0 \leq L^\prime_N - \gp^2_n \leq M(K)$.

\vskip.5cm

In order to connect $b^d_K(0)$ to the height $\gp^2_n$, it is
enough to link $b^d_K(0)$ to a seed in $T(\ell_n, \gp_n)$
and then to join it to the height $\gp^2_n$
by piling up 5 blocks of type $\bbB(t_n, \gp_{t_n})$.
As in the derivation of Theorem \ref{thm: block type 1}, the 
steering rule will be necessary to control the position
of the seeds.
This ensures that the connection from $b^d_K(0)$ to 
$T(\ell^2_n,\gp^2_n)$ occurs strictly within $\bbB(\ell^2_n,\gp^2_n)$.

Combining inequalities \eqref{eq: houba} and \eqref{eq: houba2},
we derive the Proposition with
\begin{eqnarray*}
\gep_2(\eta) = 1 - (1 -2\eta) \big( 1 - 
\eta^2 \big)^5 \leq 2\eta + 5\eta^2 \, .
\end{eqnarray*}

\end{proof}

	\subsubsection{Proof of Theorem \ref{thm: block type 2}}

Setting $\ell (K,N) = \ell^1 (K,N)$ and $h (K,N) = \gp^2 (K,N)$,
the Theorem \ref{thm: block type 2} follows by combining the FKG
inequality and Propositions \ref{pro: connections laterales} and 
\ref{pro: top connections}.


\section{Slab percolation: second renormalization step}
\label{sec: second renormalization level}

Using the renormalized blocks, we will prove that percolation occurs 
in a slab with positive probability as soon as $\bond \in \Theta_q$.
For this, we follow the cluster-growth algorithm introduced in
\cite{BGN} and prove that on a coarse grained scale large enough, 
the renormalized process dominates a supercritical two-dimensional
Bernoulli process.
This will be enough to ensure percolation in a slab.

Most of the work was done in the previous Sections to establish
uniform estimates on the occupied blocks. As a consequence, this
final step is very similar to the one devised for independent 
percolation \cite{BGN}.
With some respect, it is even simpler since the shape of the 
blocks (ratio height/width) is under control.\\

Let us fix $\bond \in \Theta_q$
and summarize the result of Theorem \ref{thm: block type 2}.
For any $\eta >0$, there exists $K$, $(L,H)$ and $(\ell,h)$ such that
\begin{eqnarray}
\label{eq: key estimate}
\FKm{{\rm f}}{\boxh^* (L,H) } 
\big(\boxh^* (\ell,h) \ \ \text{is occupied} \big) \geq 1-  
\eta  \, ,
\end{eqnarray}
where $K$ characterizes the size of the seeds and
the parameters are chosen such that
\begin{eqnarray*}
H \geq 3L, \qquad 0 \leq H - h \leq H/100  \, .
\end{eqnarray*}
At this stage the dependency on $N$ is no longer relevant. 

\medskip

We sketch below the main steps of the proof and detail only the
new features.

	\subsection{Reduction to two dimensions}

The strategy will be to create new connections by stacking rotated
translates of the block $\bbB (\ell,h)$.
In the derivation of Theorem \ref{thm: block type 1}, we already
encountered two basic features: {\bf top stacking} and {\bf steering}.
We will also present a third one: {\bf branching}.
This will enables us to define a new dependent percolation process
restricted to the slab $\slab_L = \{ -2L, \dots, 2L\}^{d-2} \times
\bbZ^2$.\\

Steering is necessary to control the deviation of the renormalized 
paths. By choosing carefully on which particular subfacet of a
block another block is attached, we may localize a trajectory.
In particular, the first $(d-2)$ coordinates happen to be irrelevant
because any stacking in a direction $\vec{e_i}$ (for $i \leq d-2$)
can be centered along the $i$-axis (we refer to \cite{BGN} for 
a complete explanation).
We stress the fact that by steering the sequence of boxes 
$\bbB (\ell,h)$ will remain inside 
$\slab_{2\ell} = \{ -2 \ell, \dots, 2 \ell \}^{d-2} \times
\bbZ^2$, nevertheless to evaluate the probability of an occupied
block one needs to average over a  bigger block $\bbB(L,H)$
(see \eqref{eq: key estimate}). 
This explain why we have to consider the renormalized process in 
the thicker slab
$\slab_{L+\ell} = \{ -L -\ell , \dots, L+\ell\}^{d-2} \times \bbZ^2$.\\

Let us now introduce some notation to emphasize the role of the
plane $(\vec{e}_{d-1},\vec{e}_d)$ which we shall refer later
on as the $(x,y)$-plane.
In order to describe the geometrical construction in the 
$(x,y)$-plane, it is convenient to rename the subfacets
of the blocks.
The direction $\vec{e}_d$ will be referred as {\bf North},
$-\vec{e}_d$ as {\bf South}, $\vec{e}_{d-1}$ as {\bf East} and
$- \vec{e}_{d-1}$ as {\bf West}.
In particular $\bbB(\ell,h)$ will be dubbed north block and denoted
by $\bbB_{\north} (\ell,h)$.
Its north face $T_{\north} (\ell,h) = T (\ell,h)$ is split into
\begin{eqnarray*}
T_{{\north},{\east}} = \{ -\ell,\ell\}^{d-2} \times \{0,\ell\} \times \{H\} \, ,\\
T_{{\north},{\west}} = \{ -\ell,\ell\}^{d-2} \times \{-\ell,0\} \times \{H\} \, .
\end{eqnarray*}
We also distinguish the western and eastern sides
\begin{eqnarray*}
S_{{\north},{\east}} = \{ -\ell,\ell\}^{d-2} \times \{\ell\} \times \{0,H\} \, ,\\
S_{{\north},{\west}} = \{ -\ell,\ell\}^{d-2} \times \{-\ell\} \times \{0,H\} \, .
\end{eqnarray*}
Rotating $\bbB_{\north}$ leads to define blocks oriented in the other 3
directions $\bbB_{\east}, \bbB_{\south}$ and $\bbB_{\west}$.
Their faces will be named by transposing the previous notation.

\medskip

We can now describe the main ingredients to concatenate paths from
different occupied blocks.
Starting with an occupied north block, {\bf top stacking} enables to 
pile up $n$ occupied blocks in the north direction with a 
probability of success at least $(1 - \eta)^n$.
As this procedure is doomed to fail sooner or later, a branching
procedure will be necessary to allow percolation in the other
directions.

\begin{figure}[h]
\begin{center}
\leavevmode
\epsfysize = 6 cm
\psfrag{0}{$b^d_K(0)$}
\psfrag{y1}[tr]{$y^{(1)}$}
\psfrag{y2}{$y^{(2)}$}
\psfrag{y3}{$y^{(3)}$}
\psfrag{z1}[Br]{$z^{(1)}$}
\psfrag{z2}[r]{$z^{(2)}$}
\psfrag{z3}{$z^{(3)}$}
\psfrag{B}[r]{$\bbB (L,H)$}
\epsfbox{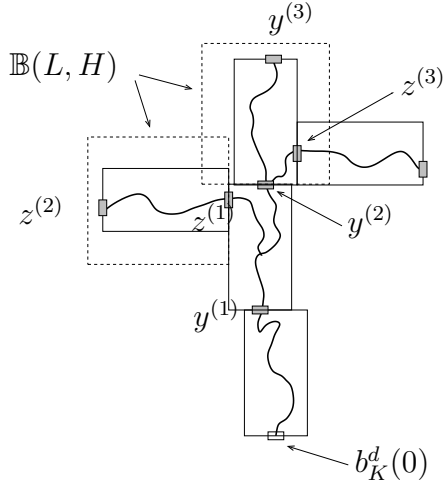}
\end{center}
\caption{A branching sequence in the North/West directions.
The domains $y^{(2)} + \bbB_{\north} (L,H)$ and $z^{(1)} + 
\bbB_{\west} (L,H)$ are depicted by the dashed blocks.}
\label{fig: branching}
\end{figure}

An occupied block contains attachment sites on each of its faces,
so that several blocks can be stacked on it. A {\bf branching} occurs
when two blocks are stacked simultaneously on the top and on one side
of an occupied block.
These blocks should be positioned in a careful way to remain
essentially independent.
For example, starting with the occupied block $\bbB_{\north} (\ell,h)$,
 a block  $y + \bbB_{\north} (\ell,h)$ can be stacked 
on a seed centered in $y \in T_{{\north},{\east}} (\ell,h)$ and a west
block  $z + \bbB_{\west} (\ell,h)$ on a seed in $S_{{\north},{\west}}$
(see figure \ref{fig: branching}).
Alternatively, exchanging East and West would lead  to a branching 
in the North/East directions.
By construction, the event that these new blocks are occupied, is
supported by the disjoint set of bonds $y + \bbB_{\north}^* (\ell,h)$
and $z + \bbB_{\west}^* (\ell,h)$.
Nevertheless to evaluate the probability of a block to be occupied
requires averaging over a larger domain $\bbB_\bullet (L,H)$
(see \eqref{eq: key estimate}).
This raises some measurability issues which will be detailed in the 
example below.

\medskip

We are going to evaluate the probability that the following sequence
of occupied blocks occurs (see figure \ref{fig: branching}).
The derivation will be similar to the one of Lemma \ref{lem: premiere connection}.
In the procedure below, the steering in the first $(d-2)$ coordinates
is implicit.

First $b^d_K (0)$ is connected strictly within $\bbB_{\north} (\ell,h)$
to a seed centered at the site $y^{(1)}$ which belongs to
 $T_{{\north},{\west}} (\ell,h)$.
The site $y^{(1)}$ is uniquely determined and is chosen as
the first site (wrt to the lexicographic order) for which the previous 
event holds.
This event is denoted by $\cE^1(y^{(1)})$. 
Thanks to the key estimate \eqref{eq: key estimate}, we have
for any domain $\gL$ containing $\bbB(L,H)$
\begin{eqnarray}
\label{eq: branch 1}
\sum_{y^{(1)}}
\FKm{{\rm f}}{\gL} \left(  \cE^1 (y^{(1)}) \right)
=
\FKm{{\rm f}}{\gL} \left(  \cC^1_K (\ell,h) \right)
\geq 
\FKm{{\rm f}}{\boxh^* (L,H) } 
\big(\boxh^* (\ell,h) \ \ \text{is occupied} \big)
\geq 1-  \eta  \, ,
\end{eqnarray}
where the sum over $y^{(1)}$ is restricted to $T_{{\north},{\west}} (\ell,h)$.

Then a branching occurs within $y^{(1)} + \bbB_{\north} (\ell,h)$, i.e.
$b^d_K (y^{(1)})$ is connected to a seed $b^d_K (y^{(2)})$ 
centered in $y^{(2)} \in  y^{(1)} + T_{{\north},{\east}} (\ell,h)$ and
 $b^d_K (y^{(1)})$ is also connected to a seed $b^{d-1}_K (z^{(1)})$
centered in $z^{(1)} \in  y^{(1)} + S_{{\north},{\west}} (\ell,h)$.
This event is denoted by $\cE^2(y^{(1)}, y^{(2)}  ,z^{(1)})$.  
As the attachment sites $y^{(2)}  ,z^{(1)}$ are uniquely determined,
the previous events are disjoint for different attachment sites.

Finally, the connections evolve in the West and North directions:
$b^{d-1}_K (z^{(1)})$ is connected to a seed $b^{d-1}_K (z^{(2)})$
in $z^{(1)} + T_{\west} (\ell,h)$ and $b^d_K (y^{(2)})$ is connected to 
a seed in $y^{(2)} + T_{{\north},{\east}} (\ell,h)$ or possibly to a 
seed in $y^{(2)} + S_{{\north},{\east}} (\ell,h)$ (on  figure
\ref{fig: branching}, a side 
stacking from $z^{(3)}$ oriented in the East direction is also depicted).
We also define
\begin{eqnarray*}
\begin{cases}
\cE^3 (y^{(2)}) = \big\{ y^{(2)}+ \bbB_{\north} (\ell,h) \ 
\text{is occupied} \big\} \, , \\ 
\cE^4 (z^{(1)}) = \big\{ z^{(1)}+ \bbB_{\west} (\ell,h)  \ 
\text{is occupied} \big\}  \, .
\end{cases}
\end{eqnarray*}

We are going to check that a branching occurs with high probability
\begin{eqnarray}
\label{eq: branching}
\sum_{y^{(1)},y^{(2)},z^{(1)}}
\FKm{{\rm f}}{\gL} \Big(  \cE^1 (y^{(1)}) \cap
\cE^2 (y^{(1)},y^{(2)},z^{(1)}) \cap \cE^3 (y^{(2)}) 
\cap \cE^4  (z^{(1)}) \Big)
\geq (1 -  \eta)^4  \, .
\end{eqnarray}
It is enough to iterate in the right order the key estimate 
\eqref{eq: key estimate}. 
Let us fix a triplet $\{ y^{(1)},y^{(2)},z^{(1)} \}$ and drop
temporarily the dependency on the sites in the events $\cE^\bullet$.
The event
$\cE^1 \cap \cE^2 \cap \cE^3$ is supported by bond configurations in 
the hyperplane $\{ x_{d-1} \geq y^{(1)}_{d-1} - \ell\}$ and
the support of $\cE^4$ lies in $\{ x_{d-1} < y^{(1)}_{d-1} - \ell \}$.
Thus conditioning outside $z^{(1)} + \bbB_{\west}(L,H)$ and using the
fact that $\cE^4 $ is non decreasing, we get
\begin{eqnarray*}
&& \FKm{{\rm f}}{\gL} \left(  \cE^1 \cap \cE^2 \cap \cE^3 \ 
\FKm{\pi}{z^{(1)} + \bbB_{\west}(L,H)} \big( \cE^4  \big) \right)\\
&& \qquad \geq \FKm{{\rm f}}{\gL} 
\left( \cE^1 \cap \cE^2 \cap \cE^3 \right)
\FKm{{\rm f}}{z^{(1)} + \bbB_{\west}(L,H)} \big( \cE^4  \big)
 \geq \FKm{{\rm f}}{\gL} 
\left( \cE^1 \cap \cE^2 \cap \cE^3 \right) \, (1 - \eta) \, .
\end{eqnarray*}
At this stage there are no more ambiguities and we estimate 
the remaining events one after the other (starting from the top)
thanks to \eqref{eq: key estimate}.
This completes \eqref{eq: branching}.

		\subsection{Cluster-growth algorithm}

We describe now the second level of renormalization.

As explained in the previous Subsection, it is enough
to consider a two-dimensional projection of the system. 
The projection of the block $\bbB_\bullet (\ell,h)$ onto the 
plane $(x,y)$ will be called a {\bf brick} and denoted
by  $\tilde \bbB_\bullet (\ell,h)$.
The orientation convention applies as well for the bricks.
To an occupied block is associated a {\bf successful} brick
with four {\bf connection sites} one in each set 
$\tilde T_{{\north},{\east}}, \tilde T_{{\north},{\west}}$,
$\tilde S_{{\north},{\east}}, \tilde S_{{\north},{\west}}$ 
which are the projections of the centers of the seeds.
A sequence of bricks build by concatenating occupied blocks 
according to the stacking procedures previously described will 
be called a {\bf successful brick sequence}.

Thus there is an immediate correspondence between blocks and
bricks and it is enough to build the second level of renormalization
in terms of bricks in $\bbZ^2$.
The second renormalization level lives on the coarse grained scale
$N = 10 L + 10 H$. 
The lattice $\bbZ^2$ is partitioned into translates of the square
$\bbS_0 = \{ - N+1 , \dots, N\}^2$ which are indexed wrt an arbitrary 
order $\{ \bbS_i \}_{i \geq 0}$. 
The algorithm is performed by inspecting the squares 
one after the other and checking iteratively some properties which 
will be detailed later on. 
If these properties are satisfied then the
 square contains a crossing cluster made of successful bricks.
A random variable $Z_i (\go)$ depending on the bond configuration
$\go$ is associated to the square $\bbS_i$.
The square $\bbS_i$ is declared to be {\bf good} if these properties
hold and we set $Z_i = 1$, otherwise $Z_i = 0$.
As it will be clear later on, any bond configuration must contain an 
open cluster intersecting all the good squares.

We explain now the construction of the renormalized process
$\{Z_i\}_i$ indexed by the squares in the domain 
$\{-M,\dots,M\}^2$ for some $M = 2mN$.
Let us suppose that $k$ squares $\bbS_0, \bbS_{i_1}, \dots, 
\bbS_{i_{k-1}}$ have been examined. Choose the next square
$\bbS_{i_k}$ as the earliest square (in the fixed ordering)
which has a face in common with one of the previous good squares.
If no such square exists then the algorithm stops, otherwise
we are going to prove that for a suitable tuning of the coarse
grained scales 
\begin{eqnarray}
\label{eq: sto domination}
\FKm{{\rm f}}{\slab_{L+\ell,M}}
\big( Z_{i_k} = 1 \, \big | \; Z_{i_{k-1}}, \dots, Z_0  \big)
\geq \ga \, ,  
\end{eqnarray}
where $\ga = \ga (\eta)$ can be chosen arbitrarily close to 1 and 
in particular much larger than the critical value of the site 
percolation in $\bbZ^2$.
We stress the fact that $\ga$ depends only on $\eta$, so that
all the parameters are determined thanks to \eqref{eq: key estimate}.
Inequality \eqref{eq: sto domination} implies that the cluster formed 
by the good squares dominates stochastically the open cluster growing 
from the origin of a supercritical site percolation process.
We refer to Grimmett, Marstrand \cite{GM} for a detailed account of 
the stochastic domination.

Before completing the construction of the renormalized process
and deriving \eqref{eq: sto domination}, we
pause to look at the proof of Theorem \ref{thm: slab threshold}.

\vskip.5cm

\noindent
{\it Proof of Theorem \ref{thm: slab threshold}.}

Fix $M = 2mN$ and choose a site $x$ in $\slab_{L+\ell, M}$. 
Let $\bbS$ be the square such that the tube 
$T_x = \{ - L-\ell, L+\ell\}^{d-2} \times \bbS$ contains $x$.
If $x =0$, then $T_0 = \{ - L-\ell, L+\ell\}^{d-2} \times \bbS_0$.
Define $\cZ(0,x)$ as the set of bond configurations such that
there exists a renormalized connected path of good squares 
from $\bbS_0$ to $\bbS$. 
By construction, one has
\begin{eqnarray}
\label{eq: renormalized path}
\FKm{{\rm f}}{\slab_{L+\ell,M}}
\big( T_0 \lra T_x   \big)
\geq \FKm{{\rm f}}{\slab_{L+\ell,M}}
\big( \cZ(0,x)  \big) \, .
\end{eqnarray}
This means that the tubes  are connected by an open path
if a connection occurs on the coarse grained level.
Thanks to the stochastic domination (see \eqref{eq: sto domination}),
one can choose $\ga = \ga(\eta)$ such that the RHS is bounded from 
below uniformly wrt to $x$ and $M$.
Therefore $\bond \geq \hat p_c (L + \ell) \geq \hat p_c$ and Theorem 
\ref{thm: slab threshold} holds.

\qed

\vskip.5cm

We turn now to the derivation of \eqref{eq: sto domination}.
As explained previously, the cluster-growth algorithm examines
in turn each square: a square $\bbS$ will be inspected only
if there exists a square $\bbS^\prime$ which shares
a common facet with $\bbS$  and
has been previously declared to be good.
In this case the state of $\bbS$ will be determined independently
of the bonds lying outside $\bbS \cup \bbS^\prime$ (more precisely,
only a small portion of $\bbS^\prime$ will be relevant).
The procedure we are going to apply is translational and rotational 
invariant, thus it is sufficient to describe it in a special case:
conditionally to the fact that $\bbS_0$ is a good square, we would 
like to determine the state of its northern neighbor 
$\bbS_1 = \bbS_0 + (0,N)$.

We define the {\bf target regions} of $\bbS_0$ as
\begin{eqnarray*}
\cT^0_{\north} = \{ -N/2, N/2\} \times \{ N - 2H \} 
\end{eqnarray*}
and the other target regions
$\cT^0_{\east},  \cT^0_{\south}, \cT^0_{\west}$ are deduced
by rotation.
In the same way, the target regions of $\bbS_1$ are obtained by
translation $\cT^1_\bullet = \cT^0_\bullet + (0,N)$.
We will see that a necessary condition for a square to be
good is to contain a successful brick 
$\tilde \bbB_\bullet (\ell,h)$ intersecting the corresponding
target region $\cT_\bullet$, where $\bullet$ ranges over the 4 
directions.\\

\begin{figure}[h]
\begin{center}
\leavevmode
\epsfysize = 8 cm
\psfrag{b0}[c]{$\tilde \bbB_0$}
\psfrag{b1}{$\tilde \bbB_1$}
\psfrag{b2}[bc]{$\tilde \bbB_2$}
\psfrag{b3}{$\tilde \bbB_3$}
\psfrag{b5}{$\tilde \bbB_5$}
\psfrag{b6}{$\tilde \bbB_6$}
\psfrag{zw}[c]{$z^{(W)}$}
\psfrag{x}{$x$}
\psfrag{y}{$y$}
\psfrag{tar}{$\cT^1_{\east}$}
\psfrag{trig}{triggering line}
\epsfbox{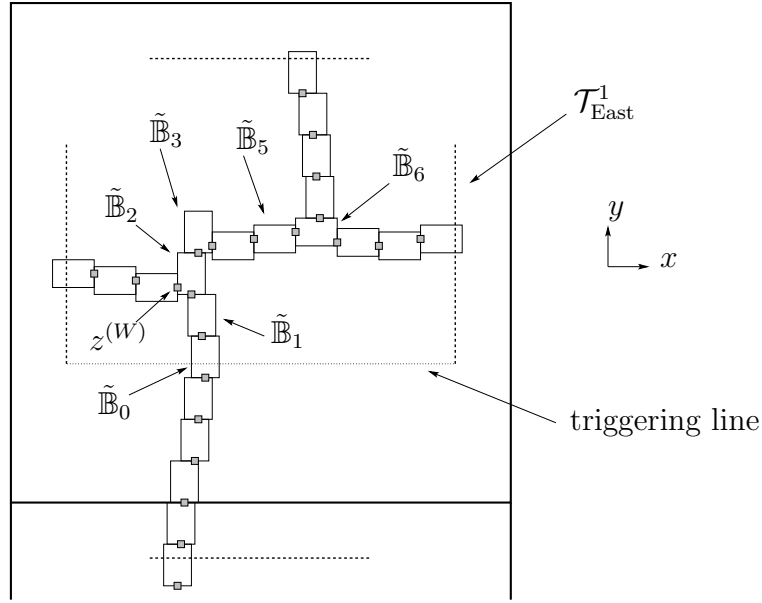}
\end{center}
\caption{A successful sequence of bricks is depicted
starting from the northern part of $\bbS_0$ and invading
$\bbS_1$. The dashed lines are the target regions.}
\label{fig: growth}
\end{figure}

Conditionally to the fact that $\bbS_0$ is a good square, 
there exists by definition a successful north brick 
$\tilde \bbB$ intersecting the target region $\cT^0_{\north}$.
This brick  will be used
as a starting point to launch connections to the target regions
of $\bbS_1$. The main ingredients are the centering
and the bifurcation rules (see Figure \ref{fig: growth}).

\vskip.5cm

\noindent
{\it Centering rules:}

First, bricks are piled up on top of $\tilde \bbB$ in the north direction 
using the steering rule in order to center the brick sequence along the
axis $x = 0$. 
At some point a brick, denoted by $\tilde \bbB_0$ will intersect the
level $\{ y = N + N/2 \}$, this triggers the bifurcation.

\vskip.5cm

\noindent
{\it Bifurcation:}

On top of $\tilde \bbB_0$ lies a connection site $z^{(1)}
= (x_1,y_1)$ (recall that $z^{(1)}$ is simply  the image of 
the center of the seed in the top face of the occupied block 
associated to $\tilde \bbB_0$).
By construction $x_1$ belongs to $\{ - 5H -5L, 5H + 5L \}$ and
we also suppose that $x_1 < 0$, the other case can be treated
by symmetry.

First a branching occurs in the North/West directions, to do so
three bricks are piled up in the north direction. 
\begin{eqnarray*}
\begin{cases}
\tilde \bbB_1 = z^{(1)} + \tilde \bbB_{\north} (\ell,h) \, , \\
\tilde \bbB_2 = z^{(2)} + \tilde \bbB_{\north} (\ell,h)
\quad \text{with} \quad
z^{(2)} \in z^{(1)} + \tilde T_{{\north},{\west}} (\ell,h) \, , \\
\tilde \bbB_3 = z^{(3)} + \tilde \bbB_{\north} (\ell,h)
\quad \text{with} \quad
z^{(3)} \in z^{(2)} + \tilde T_{{\north},{\east}} (\ell,h) \, .
\end{cases}
\end{eqnarray*}
From $\tilde \bbB_2$ a brick sequence branches in the West direction
towards the target region $\cT^1_{\west}$. We will come back to this
later.

Then an east brick $\tilde \bbB_4$ is attached to the east side of 
$\tilde \bbB_3$ and another branching occurs in the East/North
directions:
three bricks are piled up in the East direction on top of
$\tilde \bbB_4$
\begin{eqnarray*}
\begin{cases}
\tilde \bbB_4 = z^{(4)} + \tilde \bbB_{\east} (\ell,h) \quad 
\text{with} \quad z^{(4)} \in S_{{\north},{\east}} +  z^{(3)}\, , \\
\tilde \bbB_5 = z^{(5)} + \tilde \bbB_{\east} (\ell,h) \quad 
\text{with} \quad z^{(5)} \in T_{{\east},{\north}} +  z^{(4)}\, , \\
\tilde \bbB_6 = z^{(6)} + \tilde \bbB_{\east} (\ell,h) \quad 
\text{with} \quad z^{(6)} \in T_{{\east},{\north}} +  z^{(5)} \, , \\
\tilde \bbB_7 = z^{(7)} + \tilde \bbB_{\east} (\ell,h) \quad 
\text{with} \quad z^{(7)} \in T_{{\east},{\south}} +  z^{(6)} \, . 
\end{cases}
\end{eqnarray*}
From $\tilde \bbB_6$ a brick sequence branches in the North direction
towards the target region $\cT^1_{\north}$ and another sequence 
goes on in the East direction towards $\cT^1_{\east}$.

\vskip.5cm

\noindent
{\it Final connections:}

After the bifurcation, the steering rules are applied once again to
center along the axis of the target region $\cT^1_{\west}$ (resp.
North, East) the brick sequence piled up in the West direction  (resp.
North, East directions).
It remains to check that using this construction the sequences
reach the target regions and that they do not overlap.

The starting point of the western sequence $z^{(W)}= (x_W,y_W)$ 
belongs to the western side of $\tilde \bbB_3$. By construction 
\begin{eqnarray*}
x_W \in [- 2 \ell , \ell ] + x_0
, \qquad
y_W \in [ N + N/2 + H ,  N + N/2 + 3H]
\end{eqnarray*}
By piling up bricks in the western direction, the brick sequence
is localized in the strip $\{- N, x_W\} \times \{y_W -L, L\}$ and
therefore it encounters the target region $\cT^1_{\west}$ by using 
less than 20 bricks.

We proceed in the same way for the other directions and see 
that the sequences evolve in non overlapping regions.

\vskip.5cm

The inspection of a square requires less than 100 bricks,
thus the probability for $\bbS_1$ to be good  conditionally
to the fact that  $\bbS_0$ is good, is larger than 
$\ga(\eta) = (1 - \eta)^{100}$. 
Furthermore, the procedure described above does not
depend on the bond configurations outside $\bbS_0 \cap \bbS_1$.
Thus the stochastic domination inequality \eqref{eq: sto domination}
holds at any step of the cluster algorithm.

\end{document}